\documentclass[letterpaper, 11pt]{amsart}

\usepackage{amsmath,amsthm,amsfonts,amssymb,amscd}
\usepackage{bbm}
\usepackage{bm}
\usepackage{tikz}
\usepackage{tikz-cd}
\usepackage{BOONDOX-calo}
\usepackage{standalone}
\usepackage[hidelinks]{hyperref}
\usetikzlibrary{arrows,chains,matrix,positioning,scopes}
\usepackage[letterpaper, left=3cm,right=3cm, top=3cm, bottom=3cm]{geometry}

\newcommand{\CC}{{\mathbb C}}
\newcommand{\RR}{{\mathbb R}}

\newcommand{\ZZ}{{\mathbb Z}}
\newcommand{\QQ}{{\mathbb Q}}

\newcommand{\wt}{\widetilde}

\newcommand{\inEnd}{{\underline{End}}}
\newcommand{\fu}{{\mathfrak{u}}}

\newcommand{\inHom}{{\underline{Hom}}}

\newcommand\dual{\raise0.9ex\hbox{$\scriptscriptstyle\vee$}}

\theoremstyle{plain}
\newtheorem{thm}{Theorem} 
\newtheorem{prop}[thm]{Proposition}
\newtheorem{lemma}[thm]{Lemma}  
\newtheorem{cor}[thm]{Corollary}

\theoremstyle{definition}

\theoremstyle{remark}
\newtheorem{rem}[thm]{Remark}

\numberwithin{thm}{subsection}

\theoremstyle{plain}
\newtheorem{manualtheoreminner}{Theorem}
\newenvironment{thm'}[1]{%
  \renewcommand\themanualtheoreminner{#1}%
  \manualtheoreminner
}{\endmanualtheoreminner}

\tikzset{>=stealth}

\makeatletter
\def\@seccntformat#1{%
  \protect\textup{\protect\@secnumfont
    \ifnum\pdfstrcmp{subsection}{#1}=0 \bfseries\fi
    \csname the#1\endcsname
    \protect\@secnumpunct
  }%
}  
\makeatother

\newcommand\mapsfrom{\mathrel{\reflectbox{\ensuremath{\mapsto}}}}

\begin{document}

\title[The unipotent radical of the Mumford-Tate group of a very general MHS]{The unipotent radical of the Mumford-Tate group of a very general mixed Hodge structure with a fixed associated graded}
\author{Payman Eskandari and V. Kumar Murty}
\address{Department of Mathematics, University of Toronto, 40 St. George St., Room 6290, Toronto, Ontario, Canada, M5S 2E4}
\email{payman@math.toronto.edu , murty@math.toronto.edu}
\begin{abstract}
The family of all mixed Hodge structures on a given rational vector space $M_\QQ$ with a fixed weight filtration $W_\cdot$ and a fixed associated graded Hodge structure $Gr^WM$ is naturally in a one to one correspondence with a complex affine space. We study the unipotent radical of the very general Mumford-Tate group of the family. We do this by using general Tannakian results which relate the unipotent radical of the fundamental group of an object in a filtered Tannakian category to the extension classes of the object coming from the filtration. Our main result shows that if $Gr^WM$ is polarizable and satisfies some conditions, then outside a union of countably many proper Zariski closed subsets of the parametrizing affine space, the unipotent radical of the Mumford-Tate group of the objects in the family is equal to the unipotent radical of the parabolic subgroup of $GL(M_\QQ)$ associated to the weight filtration on $M_\QQ$ (in other words, outside a union of countably many proper Zariski closed sets the unipotent radical of the Mumford-Tate group is as large as one may hope for it to be). Note that here $Gr^WM$ itself may have a small Mumford-Tate group.
\end{abstract}
\maketitle

\section{Introduction}

\subsection{Statement of the results}

Let $M=(M_\QQ,W_\cdot, F^\cdot)$ be a (rational) mixed Hodge structure on the rational vector space $M_\QQ$, with weight filtration $W_\cdot$ and Hodge filtration $F^\cdot$. The Mumford-Tate group of $M$, which we denote by $G(M)$, is an algebraic subgroup of $GL(M_\QQ)$ defined more classically as the subgroup of $GL(M_\QQ)$ which fixes all the Hodge classes of weight zero (i.e. rational elements of $W_0\cap F^0$) in direct sums of the objects 
\[
M^{\otimes a,b} \ := \ M^{\otimes a}\otimes M\dual^{\otimes b} \hspace{.5in}(a,b\in \ZZ_{\geq 0}).
\]
An equivalent but perhaps more conceptual definition can be given in the language of Tannakian formalism: $G(M)$ is the fundamental group of the Tannakian category $\langle M\rangle ^{\otimes}$ of mixed Hodge structures generated by $M$, with the fiber functor being taken to be the forgetful functor. Explicitly, this means that for every object $X$ of the category generated by $M$, there is a canonical corresponding representation $\rho_X$ of the algebraic group $G(M)$ on the underlying rational vector space $X_\QQ$ of $X$, such that the forgetful functor $X\mapsto X_\QQ$ gives an equivalence of categories between $\langle M\rangle ^{\otimes}$ and the category of finite-dimensional representations of the algebraic group $G(M)$. Since $\langle M\rangle ^{\otimes}$ is generated by $M$, the map $\rho_M:G(M)\longrightarrow GL(M_\QQ)$ is injective. Identifying $G(M)$ as a subgroup of $GL(M_\QQ)$ via $\rho_M$ we obtain the same group as in the first definition. (See \cite[Lemma 1]{An92}.)
 
Since each $W_\cdot M$ is a subobject of $M$, the action of $G(M)$ respects the weight filtration on $M_\QQ$, so that $G(M)$ is contained in the parabolic subgroup $G_0(M_\QQ)$ of $GL(M_\QQ)$ associated to the weight filtration. Let $U(M)$ be the subgroup of $G(M)$ which acts trivially on $Gr^WM_\QQ$. Then $U(M)$ is the intersection of $G(M)$ with the unipotent radical $U_0(M_\QQ)$ of $G_0(M_\QQ)$; the latter group $U_0(M_\QQ)$ being the subgroup of $G_0(M_\QQ)$ which acts trivially on $Gr^WM_\QQ$. We have a split short exact sequence
\[
1 \ \longrightarrow \ U(M)  \ \longrightarrow \ G(M)  \ \longrightarrow \ G(Gr^WM)  \ \longrightarrow \ 0.  
\]
To study $G(M)$ one may then study $G(Gr^WM)$ and $U(M)$. 

Much of the interest in mixed Hodge structures is in motivic mixed Hodge structures, i.e. objects with ``geometric origin", such as the cohomology of a (not necessarily smooth or projective) complex variety. Assuming $M$ is motivic, $Gr^WM$ will be polarizable and hence $G(Gr^WM)$ reductive, so that $U(M)$ will be the unipotent radical of $G(M)$. The study of $G(Gr^WM)$ in the motivic case amounts to understanding Hodge classes in tensor powers of the pure motives $Gr^W_nM$. According to the Hodge conjecture, this is related to questions about algebraic cycles on smooth projective varieties. On the other hand, the unipotent radical $U(M)$ is more accessible: as we shall recall in the next subsection, $U(M)$ can be completely described in terms of the (Yoneda) extension classes $\mathcal{E}_p(M)$ 
\begin{equation}\label{eq8}
0 \ \longrightarrow \ W_pM  \ \longrightarrow \ M  \ \longrightarrow \ M/W_pM  \ \longrightarrow \ 0.
\end{equation}

The goal of this paper is to study the unipotent group $U(M)$ as $M$ varies through all mixed Hodge structures on the same rational vector space, with the same weight filtration and the same associated graded Hodge structure. More precisely, consider a triple of the form 
\begin{equation}\label{eq10}
\mu \ := \ (M_\QQ,W_\cdot, \widetilde{M}),
\end{equation}
where $M_\QQ$ is a fixed nonzero finite-dimensional rational vector space, $W_\cdot$ is a fixed finite increasing filtration on $M_\QQ$, and 
\[\widetilde{M}=\bigoplus\limits_{n\in\ZZ} \widetilde{M}_n,\] 
where $\widetilde{M}_n$ is a fixed pure Hodge structure of weight $n$ on the vector space $Gr^W_nM_\QQ$. Say a mixed Hodge structure $M$ is associated to $\mu$ if the underlying rational vector space of $M$ is $M_\QQ$, its weight filtration is $W_\cdot$, and 
\[
Gr^WM \ = \ \widetilde{M}.
\]
Thanks to the fact that we have fixed $Gr^WM$ in the family, the set of all mixed Hodge structures associated to $\mu$ can be identified with the set of complex points of an affine complex variety $S(\mu)$. The variety $S(\mu)$ is in fact an affine space, i.e. is isomorphic to some $\mathbb{A}_\CC^n$. To be precise, $S(\mu)$ is isomorphic as a complex variety to the quotient of a unipotent complex algebraic group by a subgroup. Such quotients are affine spaces by a theorem of Rosenlicht. (See Section \S \ref{sec: parametrizing spaces} for more details about $S(\mu)$.)

For any $M$ associated to $\mu$, the group $U(M)$ is contained in $U_0(M_\QQ)$ (the subgroup of $GL(M_\QQ)$ which respects the weight filtration and acts trivially on $Gr^W(M_\QQ)$). We have equality of the two groups if and only if their Lie algebras coincide, i.e.
\[
\fu(M) \ := \ Lie(U(M) \ = \ W_{-1}End(M_\QQ) \ := \ \{f\in End(M_\QQ) \, :\, f(W_{n}M_\QQ)\subset W_{n-1}M_\QQ \ (\forall n)\}
\] 
($W_{-1}End(M_\QQ)$ being the Lie algebra of $U_0(M_\QQ)$). Following \cite{EM21b}, let us say $\fu(M)$ is {\it large} if it is equal to $W_{-1}End(M_\QQ)$. Let $D(\mu)$ be the subset of $S(\mu)(\CC)$ consisting of all those $M$ whose $\fu(M)$ is not large. A main result of the paper is the following: 

\begin{thm'}{A}\label{main thm}
Let $\mu$ be as above. Suppose the following two conditions hold:
\begin{itemize}
\item[(i)] $\widetilde{M}$ is semisimple.
\item[(ii)] For any two distinct pairs $(i,j)$ and $(i',j')$ of integers with $i<j$ and $i'<j'$, the two objects
\[
\inHom(\widetilde{M}_j,\widetilde{M}_i) \quad \text{and} \quad \inHom(\widetilde{M}_{j'},\widetilde{M}_{i'}) 
\] 
have no nonzero isomorphic subobjects. ($\inHom$ means internal Hom in the category of mixed Hodge structures.) 
\end{itemize}
Then $D(\mu)$ is contained in a union of countably many proper Zariski closed subsets of $S(\mu)(\CC)$.
\end{thm'}

In other words, assuming Conditions (i) and (ii), for a very general $M$ associated to $\mu$ the unipotent radical of the Mumford-Tate group of $M$ is large. Note that Condition (ii) is satisfied for example if the numbers 
\[
i-j \hspace{.5in}(i,j\in\ZZ,\,  i<j, \, \wt{M}_i\neq 0, \, \wt{M}_j\neq 0)
\]
are all distinct.

In the process of proving Theorem \ref{main thm} we shall show that the jumping Hodge locus of the family (where one has ``extra" Hodge tensors) is a union of countably many Zariski closed subsets of $S(\mu)(\CC)$ (see \S \ref{sec: method of proofs in Intro} for more details). But the result says more, namely that outside this jumping locus the unipotent radical of the Mumford-Tate group is equal to $U_0(M_\QQ)$. Note that here (and in all the results below) the Mumford-Tate group of the fixed associated graded $Gr^WM$ ( = $\widetilde{M}$) may be small. A fact to keep in mind in contextualizing these statements is that the universal family of mixed Hodge structures associated to $\mu$ is not necessarily a variation of mixed Hodge structures, as it might fail Griffiths' transversality. 

Let $p$ be an integer. For any mixed Hodge structure $M$ in our family, let
\[
\fu_p(M) \ = \ \fu(M)\cap Hom(M_\QQ/W_pM_\QQ, W_pM_\QQ),
\]
where $Hom(M_\QQ/W_pM_\QQ, W_pM_\QQ)$ is considered as a subspace of $W_{-1}End(M_\QQ)$ in the obvious way, using functoriality properties of of $Hom$. Let us say $\fu_p(M)$ is large if it is equal to $Hom(M_\QQ/W_pM_\QQ, W_pM_\QQ)$. Then $\fu(M)$ is large if and only if $\fu_p(M)$ is large for all $p$. Let $D^p(\mu)$ be the subset of $S(\mu)(\CC)$ consisting of all those $M$ whose $\fu_p(M)$ is not large. We shall prove the following result, from which Theorem \ref{main thm} follows immediately:

\begin{thm'}{B}\label{thm2}
Let $\mu$ be as above. Fix $p\in\ZZ$. Suppose the following two conditions hold:
\begin{itemize}
\item[(i)] $\widetilde{M}$ is semisimple.
\item[(ii)] The two objects 
\[
\bigoplus\limits_{(i,j)\in J_1^p}\inHom(\widetilde{M}_j,\widetilde{M}_i)\quad \text{and}\quad \bigoplus\limits_{(i,j)\in J_2^p}\inHom(\widetilde{M}_j,\widetilde{M}_i)
\]
where 
\begin{align*}
J_1^p \ &:= \ \{(i,j)\in \ZZ^2 \, : \, i\leq p<j\} \\
J_2^p \ &:= \ \{(i,j)\in \ZZ^2 \, : \, i<j\} \, \setminus  J_1^p
\end{align*}
do not have any nonzero isomorphic subobjects.
\end{itemize}
Then $D^p(\mu)$ is contained in a union of countably many proper Zariski closed subsets of $S(\mu)(\CC)$. 
\end{thm'}

Condition (ii) in both results above is an ``independence axiom" requirement, in the sense of \cite{EM21b}. We shall explain the relevance of Conditions (i) and (ii) to the proof of Theorem \ref{thm2} (and hence Theorem \ref{main thm}) in \S \ref{sec: method of proofs in Intro} below. Most likely, the assertions made in the theorems remain true even in the absence of these conditions, although our proof needs them (more on this is discussed in \S \ref{sec: method of proofs in Intro}). We will however prove a simpler result in the same spirit without assuming Conditions (i) or (ii), which we now explain.

Continue to fix $p$. The filtration $W_\cdot$ on $M_\QQ$ induces filtrations on $W_pM_\QQ$ and $M_\QQ/W_pM_\QQ$, which with abuse of notation we also denote by $W_\cdot$. Then $W_p\widetilde{M}$ and $\widetilde{M}/W_p\widetilde{M}$  are graded Hodge structures on $W_pM_\QQ$ and $M_\QQ/W_pM_\QQ$, respectively. Set
\begin{align*}
\mu_p \ &:= \ (W_pM_\QQ, W_\cdot, W_p\widetilde{M}) \\
\mu_{>p} \ &:= \ (M/W_pM_\QQ, W_\cdot, \widetilde{M}/W_p\widetilde{M}).  
\end{align*}
These are triples of a similar type to $\mu$. There is a morphism 
\begin{equation}\label{eq32}
\Theta_p\  : \ S(\mu) \ \longrightarrow S(\mu_p) \times S(\mu_{>p})
\end{equation}
which at the complex level sends a mixed Hodge structure $M$ on $M_\QQ$ associated to $\mu$ to the pair of mixed Hodge structures $(W_pM, M/W_pM)$ attached to $\mu_p$ and $\mu_{>p}$. Moreover, each of its fibers is seen to be an affine space, of positive dimension as long as $W_pM_\QQ$ and $M_\QQ/W_pM_\QQ$ are nonzero. Given $(x,y)\in S(\mu_p)(\CC) \times S(\mu_{>p})(\CC)$, let $S(\mu)(\CC)_{(x,y)}$ be the set of complex points of the fiber above $(x,y)$ (parametrizing mixed Hodge structures $M$ associated to $\mu$ whose $W_pM$ and $M/W_pM$ respectively are the mixed Hodge structures associated to $\mu_p$ and $\mu_{>p}$ corresponding to the points $x$ and $y$). Let $D^p(\mu)_{(x,y)}$ be the subset of $S(\mu)(\CC)_{(x,y)}$ consisting of those $M$ in $S(\mu)(\CC)_{(x,y)}$ whose $u_p(M)$ is not large. We show that:

\begin{thm'}{C}\label{thm1}
Fix $p\in\ZZ$ and
\[
(x,y) \ \in \ S(\mu_p)(\CC) \times S(\mu_{>p})(\CC).
\]
Then $D^p(\mu)_{(x,y)}$ is a union of countably many proper Zariski closed subsets of $S(\mu)(\CC)_{(x,y)}$. 
\end{thm'}

\subsection{The method of the proofs}\label{sec: method of proofs in Intro}
The more classical approach to the study of Mumford-Tate groups is through studying their invariants (e.g. see \cite{De82}). Here, since we are interested in the unipotent radical only, instead we exploit the relations between the unipotent group $U(M)$ and the extensions $\mathcal{E}_p(M)$ of Eq. \eqref{eq8}, which we shall briefly recall here. The rational vector space $\fu(M)$ (resp. $\fu_p(M)$) underlies a mixed Hodge substructure of $W_{-1}\inEnd(M)$ (resp. $\inHom(M/W_pM,W_pM)$), which we denote by $\underline{\fu}(M)$ (resp. $\underline{\fu}_p(M)$). By a result of Deligne from \cite[Appendix]{Jo14}, $\underline{\fu}(M)$ is determined by the extension
\[
\mathcal{E}(M) \ := \ \sum\limits_p \mathcal{E}_p(M) \ \in \ Ext(\mathbbm{1}, W_{-1}\inEnd(M)),
\]
where each $\mathcal{E}_p(M)$ is first considered as an extension of $\mathbbm{1}$ by $W_{-1}\inEnd(M)$ through pushforward along the inclusion $\inHom(M/W_pM,W_pM)\subset W_{-1}\inEnd(M)$. Indeed, Deligne proves that $\underline{\fu}(M)$ is the smallest subobject of $W_{-1}\inEnd(M)$ such that the pushforward of $\mathcal{E}(M)$ along
\[
W_{-1}\inEnd(M) \ \longrightarrow W_{-1}\inEnd(M)/\underline{\fu}(M)
\] 
splits. By weight considerations, splitting of this pushforward can be translated as the existence of a nonzero Hodge class of weight zero in the quotient
\[
W_{0}\inEnd(M)^\dagger/\underline{\fu}(M),
\]
where $W_{0}\inEnd(M)^\dagger$ is a certain subobject of $W_{0}\inEnd(M)$, sitting as the middle object of the extension $\mathcal{E}(M)$. On the other hand, by \cite[Theorem 3.3.1]{EM21a} (which is obtained by a small modification of the proof of Hardouin's \cite[Theorem 2]{Har11} and \cite[Theorem 2.1]{Har06}, those being generalizations of a result of Bertrand from \cite{Be01}), the extension $\mathcal{E}_p(M)$ by itself determines $\underline{\fu}_p(M)$, although in a more complicated way: Use canonical constructions to write $\mathcal{E}_p(M)$ as the extension
\[
0 \ \longrightarrow \ \inHom(M/W_pM,W_pM)  \ \longrightarrow \ \inHom(M/W_pM,M)^\dagger  \ \longrightarrow \ \mathbbm{1}  \ \longrightarrow \ 0.
\]
The object $\inHom(M/W_pM,M)^\dagger$ sitting in the middle is the subobject of $\inHom(M/W_pM,M)$ whose underlying rational vector space denoted by $Hom(M_\QQ/W_pM_\QQ, M_\QQ)^\dagger$ consists of all linear maps $f: M_\QQ/W_pM_\QQ\longrightarrow M_\QQ$ which after composing with the quotient $M_\QQ\longrightarrow M_\QQ/W_pM_\QQ$ are a scalar multiple $\lambda(f)$ of the identity on $M_\QQ/W_pM_\QQ$. The surjective arrow in the sequence is $f\mapsto \lambda(f)$. Then $\underline{\fu}_p(M)$ is the smallest subobject of $\inHom(M/W_pM,W_pM)$ such that
\begin{equation}\label{eq9}
\inHom(M/W_pM,M)^\dagger/\underline{\fu}_p(M)
\end{equation}
belongs to the Tannakian subcategory of the category of mixed Hodge structures generated by $M/W_pM$ and $W_pM$. 

This last result on the relation between $\mathcal{E}_p(M)$ and $\underline{\fu}_p(M)$ enables one to prove Theorem \ref{thm1} in a fairly straightforward way, the key being that in each fiber $S(\mu)(\CC)_{(x,y)}$ both $M/W_pM$ and $W_pM$ are fixed, so that up to isomorphism there are only countably many possibilities for $\inHom(M/W_pM,M)^\dagger/\underline{\fu}_p(M)$.

The relationship between $\mathcal{E}_p(M)$ and $\underline{\fu}_p(M)$ simplifies in the presence of Conditions (i) and (ii) of Theorem \ref{thm2}. Indeed, by \cite[Corollary 5.3.2]{EM21b}, assuming Conditions (i) and (ii), $\underline{\fu}_p(M)$ is the smallest subobject of $\inHom(M/W_pM,W_pM)$ such that the pushforward of $\mathcal{E}_p(M)$ along the quotient
\[
 \inHom(M/W_pM,W_pM) \ \longrightarrow \  \inHom(M/W_pM,W_pM)/\underline{\fu}_p(M)
\]
splits. By weight considertaions, this splitting is equivalent to existence of a nonzero Hodge class of weight zero in Eq. \eqref{eq9}. 

To prove Theorem \ref{thm2}, we use the previous result from \cite{EM21b}, together with some geometric results about algebraic families of mixed Hodge structures associated to a triple of the form Eq. \ref{eq10}. By an algebraic family of mixed Hodge structures associated to a triple $\mu$ as above parametrized by a complex variety $X$ we mean a morphism 
\[
 j \, : \, X \ \longrightarrow \ S(\mu).
\]
Some examples of algebraic families are follows:
\begin{itemize}
\item[-] Taking $X=S(\mu)$ and $j$ the identity map we obtain an algebraic family associated to the triple $\mu$; we call this family the {\it universal} family. 
\item[-] Families defined using canonical linear algebraic constructions involving $Hom$, directs sums, tensors and duals, and the weight filtration: Given the triple $\mu$ of Eq. \eqref{eq10}, it is easy to see that $\{\inHom(M/W_pM, W_pM)\}_{M\in S(\mu)(\CC)}$ is an algebraic family parametrized by $S(\mu)$, associated to the triple
\[
(Hom(M_\QQ/W_pM_\QQ, W_pM_\QQ), W_\cdot, \inHom(\widetilde{M}/W_p\widetilde{M},W_p\widetilde{M})).
\]
Here, as well as below, the weight filtration on $Hom(M_\QQ/W_pM_\QQ, W_pM_\QQ)$ or other objects involving $M_\QQ$ is defined in the natural way using the weight filtration on $M_\QQ$. Similarly, we have algebraic families
\[
(\inHom(M/W_pM, M)^\dagger)_{M\in S(\mu)(\CC)}
\]
and 
\[
(W_{-1}\inEnd(M))_{M\in S(\mu)(\CC)}
\]
parametrized by $S(\mu)$ associated respectively to 
\[
(Hom(M_\QQ/W_pM_\QQ, M_\QQ)^\dagger, W_\cdot, \inHom(\widetilde{M}/W_p\widetilde{M},\widetilde{M})^\dagger)
\]
and
\[
(W_{-1}End(M_\QQ), W_\cdot, W_{-1}\inEnd(\widetilde{M})).
\]
More generally, we can similarly get families parametrized by $S(\mu)$ associated to various triples that can be produced from $\mu$ using operations of tensor products, $Hom$, dualizing, direct sums, applying and modding out by $W_p$. 
\end{itemize}

An algebraic family of mixed Hodge structures associated to a triple (in the above sense) parametrized by $X$
naturally gives rise to a family of mixed Hodge structures on $X$ (or rather $X(\CC)$) in the usual sense (i.e. a holomorphic vector bundle with a flat connection, filtered by weight and Hodge filtrations with the weight filtration being flat with respect to the connection, and such that the fibers equipped with the induced filtrations are mixed Hodge structures). However, it is important to note that the families of interest to us in the paper, in particular the families of the examples above, are usually not variations of mixed Hodge structures (that is, they may not satisfy the Griffiths transversality condition). 

In order to prove Theorem \ref{thm2} we need a closed Hodge loci result for our algebraic families. The theorem of Cattani-Deligne-Kaplan \cite{CDK} does not apply to our setting, since as our families are usually not variations. But again, thanks to the fact that $Gr^WM$ is fixed in the family, an analogous statement can be established as a consequence of Kostant-Rosenlicht's closed orbit theorem for actions of unipotent groups on affine varieties. To be more precise, given an algebraic family of mixed Hodge structures associated to $\nu=(N_\QQ,W_\cdot, \widetilde{N})$ and parametrized by $X$, and given a subspace $A_\QQ\subset N_\QQ$, let $X(\CC)_{A_\QQ}$ be the set of points in $X(\CC)$ consisting of $x\in X(\CC)$ such that $A_\QQ$ underlies a subobject of $N_x$, where $N_x$ is the mixed Hodge structure associated to $\nu$ at $x$. As an application of Kostant-Rosenlicht's closed orbit theorem we 
we show that $X(\CC)_{A_\QQ}$ is a Zariski closed subset of $X(\CC)$. In particular, given any $v\in N_\QQ$, the Hodge locus\footnote{Note that we use the term {\it Hodge locus} in a different way (although related, of course) to the more common usage of the term.} 
\[
X(\CC)_v \ := \ \{x\in X(\CC) \, : \, \text{$v$ is a Hodge class of weight zero in $N_x$}\}
\]
is a Zariski closed subset of $X(\CC)$. 

Given an algebraic family of mixed Hodge structures associated to $\nu=(N_\QQ,W_\cdot, \widetilde{N})$ and parametrized by $X$, let us say a subspace $A_\QQ\subset N_\QQ$  is a global Hodge subspace for the family
if $X(\CC)_{A_\QQ}=X(\CC)$. If $A_\QQ$ is a global Hodge subspace, denote the mixed Hodge structure on $A_\QQ$ at $x\in X(\CC)$ by $A_x$, so that $A_x$ is the subobject of $N_x$ with underlying rational space $A_\QQ$. The final ingredient of the proof of Theorem \ref{thm2} is that when $A_\QQ\subset N_\QQ$ is a global Hodge subspace, the mixed Hodge structures
\[
(N_x/A_x)_{x\in X(\CC)}
\] 
form an algebraic family of mixed Hodge structures associated to the triple
\begin{equation}\label{eq11}
(N_\QQ/A_\QQ, W_\cdot, \widetilde{N}/\widetilde{A}),
\end{equation}
where $\widetilde{A}$ is the subobject of $\widetilde{N}$ with underlying rational vector space $Gr^WA_\QQ$. The nontrivial part of this assertion is the claim about algebraicity. This is proved using the fact that in characteristic zero any torsor of a unipotent algebraic group over an affine variety is trivial. 

With these in hand, the proof of Theorem \ref{thm2} goes as follows. Given $\mu$ as in Eq. \eqref{eq10}, it suffices to show that for every proper subspace $A_\QQ$ of $Hom(M_\QQ/W_pM_\QQ,W_pM_\QQ)$, the set of all those $s$ in $S(\mu)(\CC)$ such that $A_\QQ$ underlies a subobject of $\inHom(M_s/W_pM_s,W_pM_s)$ and $\fu_p(M_s)\subset A_\QQ$ is contained in a countable union of proper Zariski closed subsets of $S(\mu)(\CC)$. It suffices to restrict our attention to those $A_\QQ$ that are global Hodge subspaces in the family $(\inHom(M_s/W_pM_s,M_s)^\dagger)_{s\in S(\mu)(\CC)}$ (as otherwise, $S(\mu)(\CC)_{A_\QQ}$ for the family is already a proper Zariski closed subset of $S(\mu)(\CC)$).
For any such $A_\QQ$, now we have an algebraic family 
\[
(\inHom(M_s/W_pM_s,M_s)^\dagger/A_s)_{s\in S(\mu)(\CC)}
\]
parametrized by $S(\mu)$. Assuming Conditions (i) and (ii) of Theorem \ref{thm2}, by the characterization of $\underline{\fu}_p$ from \cite{EM21b} recalled above we can write the set of those $s\in S(\mu)(\CC)$ such that $\underline{\fu}_p(M_s)\subset A_s$ as a countable union of Zariski closed subsets of $S(\mu)(\CC)$: the union is over the nonzero elements of 
\[
Hom(M_\QQ/W_pM_\QQ, M_\QQ)^\dagger/A_\QQ \, ,
\]
and the corresponding closed set for each element is the Hodge locus of the element. It remains to show that each of these Hodge loci is a proper subset of $S(\mu)(\CC)$, which is not difficult, in particular because the Lie subalgebra
\[
Hom(M_\QQ/W_pM_\QQ, W_pM_\QQ)
\]
of $W_{-1}End(M_\QQ)$ is abelian.

As the reader can see from this argument, the reason we need Hypotheses (i) and (ii) in Theorem \ref{thm2} (and hence Theorem \ref{main thm}) is that our argument uses \cite[Corollary 5.3.2]{EM21b} relating 
\[\mathcal{E}_p(M) \ \in \ Ext(\mathbbm{1}, \inHom(M/W_pM,W_pM))\]
and $\underline{\fu}_p(M)$ for a mixed Hodge structure $M$, rather than Deligne's characterization of $\underline{\fu}(M)$ in terms of the total extension 
\[\mathcal{E}(M) \ \in \ Ext(\mathbbm{1}, W_{-1}\inEnd(M)).\]
If one were to repeat the same argument for the family $(W_0\inEnd(M_s)^\dagger)_{s\in S(\mu)(\CC)}$ and use Deligne's result (and hence remove Hypotheses (i) and (ii) from the results), the argument goes through until the very last step, where one needs to show that the Hodge locus of every nonzero element of
\[
W_{0}End(M_\QQ)^\dagger/A_\QQ
\]
is a {\it proper} subset of $S(\mu)(\CC)$, where $A_\QQ\subsetneq W_{-1}End(M_\QQ)$ is a global Hodge subspace for the family $(W_{0}\inEnd(M_s)^\dagger)_{s\in S(\mu)(\CC)}$. The difficulty in showing this last step is because of the nonabelian nature of $W_{-1}End(M_\QQ)$.

\subsection{Reinterpretation of Theorem \ref{main thm} in terms of invariants}
In view of the result about Zariski closedness of Hodge loci in algebraic families and the relation between Mumford-Tate groups and their invariants, we can reformulate the assertion of Theorem \ref{main thm} as follows. Given the triple $\mu$, we consider algebraic families of the form
\begin{equation}\label{eq31}
(\bigoplus\limits_{(a,b)\in I} \, M^{\otimes a,b}) \, _{M\in S(\mu)(\CC)},
\end{equation}
where $I$ is a finite subset of $\ZZ_{\geq 0}\times \ZZ_{\geq 0}$. Consider the union
\[
\bigcup\limits_{I, v} \, S(\mu)(\CC)_v,
\]
where the union is over all finite subsets $I$ of $\ZZ_{\geq 0}\times \ZZ_{\geq 0}$ and all $v \in \bigoplus_{(a,b)\in I} \, M_\QQ^{\otimes a,b}$ such that the Hodge locus $S(\mu)(\CC)_v$ is not all of $S(\mu)(\CC)$ (i.e. $v$ is not a global Hodge class for the family Eq. \eqref{eq31}). This is a countable union of proper Zariski closed subsets of $S(\mu)(\CC)$, and outside it the only Hodge classes in each of the families of the form Eq. \eqref{eq31} are the global Hodge classes. Hence the Mumford-Tate group $G(M)$ for a very general $M$ is the family $S(\mu)$ is the subgroup of $GL(M_\QQ)$ which fixes all the global Hodge classes $v$ in the families of the form Eq. \eqref{eq31}. Thus the assertion of Theorem \ref{main thm} is equivalent to the following statement: for any finite $I\subset\ZZ_{\geq 0}\times \ZZ_{\geq 0}$, if 
\[
v \ \in \ \bigoplus\limits_{(a,b)\in I} \, M_\QQ^{\otimes a,b}
\] 
is a global Hodge class for the family
\[
(\bigoplus\limits_{(a,b)\in I} \, M^{\otimes a,b})\, _{M\in S(\mu)(\CC)},
\]
then $v$ is fixed by $U_0(M_\QQ)$. 

\subsection{Organization of the paper}
The paper is organized as follows. In the next section, we shall construct $S(\mu)$ and study the map $\Theta_p$ of Eq. \eqref{eq32}. In \S \ref{sec: algebraic families attached to triples} we study algebraic families of mixed Hodge structures associated to a triple; the main statements to be proved will be the closedness of Hodge loci and the result about quotients by global Hodge subspaces. In \S \ref{sec: proof of main thm} (resp. \S \ref{sec: proof of Thm C}) we use the results of the previous sections together with \cite{EM21b} to establish Theorem \ref{thm2} (resp. \ref{thm1}).

\section{Parametrizing spaces}\label{sec: parametrizing spaces}

\subsection{Definition of $S(\mu)$}\label{sec: def of S}
Throughout the paper, we shall constantly work with triples of the form
\[
(M_\QQ, W_\cdot, \widetilde{M}),
\] 
where 
\begin{itemize}
\item[(i)] $M_\QQ$  is a finite-dimensional vector space over $\QQ$, 
\item[(ii)] $W_\cdot$ is a finite increasing filtration  on $M_\QQ$, and
\item[(iii)] $\wt{M}=\bigoplus\limits_{n\in\ZZ} \wt{M}_n$, where for each $n$, $\wt{M}_n$ is a Hodge structure of weight $n$ on the rational vector space $Gr^W_nM_\QQ$.  
\end{itemize}
We will simply call such data a {\it triple}. We will use lower case Greek letters $\mu$ and $\nu$ to name our triples.

Given a triple $\mu=(M_\QQ, W_\cdot, \wt{M})$, we define the following objects:
\begin{itemize}
\item $U(\mu)$: This is the unipotent radical of the parabolic subgroup of $GL(\widetilde{M}_\CC)$ associated to the weight filtration on $\widetilde{M}_\CC$; thus $U(\mu)$ is a complex algebraic group, the subgroup of $GL(\widetilde{M}_\CC)$ whose group of $R$-valued points for any commutative $\CC$-algebra $R$ is the group of automorphisms $\alpha$ of  $\widetilde{M}_\CC\otimes R$ that preserve the weight filtration, and moreover induce the identity map on $Gr^W\widetilde{M}_\CC$. Thus writing the elements of $GL(\widetilde{M}_\CC)$ as block matrices corresponding to the decomposition $\widetilde{M}_\CC=\bigoplus_n (\widetilde{M}_n)_\CC$, the group $U(\mu)$ is just the subgroup of block upper triangular matrices with identity maps on the diagonal. The Lie algebra of $U(\mu)$ is
\[
W_{-1}End(\widetilde{M}_\CC) \ := \ \{f\in End(\widetilde{M}_\CC) \, : \, \forall n\in\ZZ, \, f(W_n\widetilde{M}_\CC)\subset W_{n-1}\widetilde{M}_\CC\}.
\]
(If we write the endomorphisms as block matrices this is the space of strictly upper triangular elements of $End(\widetilde{M}_\CC)$. Note that $W_{-1}End(\widetilde{M}_\CC)$ is the underlying rational vector space of $W_{-1}\inEnd(\wt{M})$.)
\item $F^0U(\mu)$: The subgroup of $U(\mu)$ consisting of the elements which preserve the Hodge filtration on $\widetilde{M}_\CC$ (for $R$-valued points, they preserve the Hodge filtration after extending the scalars to $R$). The Lie algebra of $F^0U(\mu)$ is 
\[
F^0W_{-1}End(\widetilde{M}_\CC) \ := \ \{f\in W_{-1}End(\widetilde{M}_\CC) \, : \, \forall p\in\ZZ, \, f(F^p\widetilde{M}_\CC)\subset F^p\widetilde{M}_\CC\}.
\]
Note that $F^0W_{-1}End(\widetilde{M}_\CC)$ is $F^0$ of the Hodge structure $W_{-1}\inEnd(\widetilde{M})$. In particular, since $W_{-1}\inEnd(\widetilde{M})$ has weights $\leq -1$, we have $U(\mu)=F^0U(\mu)$ if and only if $W_{-1}\inEnd(\widetilde{M})$ is zero, i.e. if and only if $\widetilde{M}$ is pure, or in other words, either $M_\QQ$ is zero or its weight filtration has only one step. (In those situations, every thing we will be doing will be trivial.)
\item $T(\mu)$: This is a complex variety with the following functor of points: for any commutative $\CC$-algebra $R$, the set of $R$-valued points of $T(\mu)$ is the set of all $R$-linear isomorphisms 
\[
\alpha \, : \, \widetilde{M}_\CC\otimes R \ \longrightarrow \ M_\CC\otimes R
\]
which satisfy the following two conditions:
\begin{itemize}
\item[(i)] $\alpha$ preserves the weight filtration.
\item[(ii)] Identifying $Gr^W\widetilde{M}_\CC\otimes R$ with $\widetilde{M}_\CC\otimes R$ via the canonical isomorphism (since $\widetilde{M}_\CC\otimes R$ is graded by weight), the induced map 
\[\begin{array}{ccc}
Gr^W\widetilde{M}_\CC\otimes R & \ \xrightarrow{Gr^W\alpha} \ & Gr^WM_\CC\otimes R\\[3pt]
\rotatebox{90}{$=$} & & \rotatebox{90}{$=$}\\[3pt]
\widetilde{M}_\CC\otimes R & & \widetilde{M}_\CC\otimes R
\end{array}
\]
is the identity map.
\end{itemize}
Note that $T(\mu)$ is a right $U(\mu)$-torsor (the action given by precomposition). In particular, $T(\mu)$ is an affine space.
\item Finally, set 
\[
S(\mu) \ := \ T(\mu)/F^0U(\mu).
\]
This is (non-canonically) isomorphic to $U(\mu)/F^0U(\mu)$, and hence being isomorphic to a quotient of a unipotent algebraic group by a subgroup it is an affine space (by \cite[Theorem 5]{Ro63}). Moreover, as long as the weight filtration has more than one step, $U(\mu)/F^0U(\mu)$ will have dimension $\geq 1$.
\end{itemize}

\begin{rem}\label{elements of T and sections of W_n -> Gr_n}
Let $R$ be a commutative $\CC$-algebra . Given an $R$-linear map $\alpha: \wt{M}_\CC\otimes R\longrightarrow M_\CC\otimes R$, write $\alpha=\sum\limits_n \alpha_n$, where $\alpha_n$ is the restriction of $\alpha$ to the weight $n$ component of $\wt{M}_\CC\otimes R$, i.e. $(\wt{M}_n)_\CC\otimes R$. Then $\alpha$ is in $T(\mu)(A)$ if and only if each $\alpha_n$ is a section of the natural map
\begin{equation}\label{eq17}
W_nM_\CC\otimes R \ \longrightarrow \ (\wt{M}_n)_\CC\otimes R \, = \, Gr^W_nM_\CC\otimes R.
\end{equation}
Here and in what follows, with abuse of notation a map into $W_nM_\CC\otimes R$ and its composition with the inclusion $W_nM_\CC\otimes R\subset M_\CC\otimes R$ are denoted by the same symbol. We get a bijection
\begin{align}
T(\mu)(R) \ &\cong \ \prod_n \, Sec\, (W_nM_\CC\otimes R \ \longrightarrow \ (\wt{M}_n)_\CC\otimes R)  \label{eq19}\\
\alpha  \ & \mapsto \ (\alpha_n)_n \notag \\
\sum_n\, \alpha_n \ & \mapsfrom \ \ (\alpha_n)_n,\notag
\end{align}
where $Sec \, (\pi)$ means the set of (linear) sections of a surjective linear map $\pi$. We get an isomorphism of complex varieties
\[
T(\mu) \ \cong \  \prod_n \, Sec \, (W_nM_\CC \ \longrightarrow \ (\wt{M}_n)_\CC)\, ^a,
\]
where $Sec \, (W_nM_\CC \ \longrightarrow \ (\wt{M}_n)_\CC)\, ^a$ means the $\CC$-scheme associated to the sections of the map (with the set of $R$-valued points equal to the set of sections of $W_nM_\CC\otimes R \ \longrightarrow \ (\wt{M}_n)_\CC\otimes R$).

It is easy to see that for every $\alpha\in T(\mu)(R)$ and every $n$, we have
\[
W_nM_\CC\otimes R \ = \ \sum\limits_{m\leq n} Im(\alpha_m),
\]
where the sum on the right is direct.
\end{rem}

\subsection{Mixed Hodge structures associated to a triple}\label{sec: S parametrizes the MHSs associated to the triple}
Let $\mu=(M_\QQ, W_\cdot, \wt{M})$ be a triple. We say a mixed Hodge structure $N$ is associated to $\mu$ if the underlying rational vector space of $N$ is $M_\QQ$, its weight filtration is $W_\cdot$ given in the triple, and 
\[
Gr^WN \ = \ \wt{M}.
\]
Note that any mixed Hodge structure associated to $\mu$ is determined by the corresponding Hodge filtration on $M_\CC$. The requirement of being a mixed Hodge structure associated to $\mu$ translates in terms of the corresponding filtration on $M_\CC$ to the following: the filtration on $M_\CC$ induces the fixed Hodge filtration on $\wt{M}_\CC$ by the formula 
\[F^\cdot (\wt{M}_n)_\CC \ = \ Im\bigm(F^\cdot W_nM_\CC \, \hookrightarrow \, W_nM_\CC \, \longrightarrow \, W_nM_\CC/W_{n-1}M_\CC \, =(\wt{M}_n)_\CC\bigm)\]
for all $n$. Here, in its first appearance, $F^\cdot$ is the Hodge filtration on $(\wt{M}_n)_\CC$ corresponding to the pure Hodge structure $\wt{M}_n$ (which is a part of the defining data of $\mu$), and in its second appearance $F^\cdot$ is the filtration on $M_\CC$. (Note that if a filtration $F^\cdot$ on $M_\CC$ satisfies the formula, then $(M_\QQ,W_\cdot, F^\cdot)$ is a mixed Hodge structure, since by assumption the induced filtrations on the $Gr^W_n M_\CC$ make $Gr^W_n M_\QQ$ into a Hodge structure of weight $n$.)

The collection of all mixed Hodge structures associated to $\mu$ is in a canonical bijection with the set of complex points of $S(\mu)$, as we shall see below. Given any $\alpha\in T(\mu)(\CC)$, we may use $\alpha$ to transport the Hodge filtration of $\widetilde{M}_\CC$ to $M_\CC$: set
\begin{equation}\label{eq18}
F^\cdot_\alpha M_\CC \ := \ \alpha(F^\cdot \widetilde{M}_\CC) \, = \, \sum\limits_n \alpha_n(F^\cdot (\widetilde{M}_n)_\CC)\, , 
\end{equation}
where the sum on the right is direct. Let
\[
M_\alpha \ := \ (M_\QQ, W_\cdot, F^\cdot_\alpha).
\]

\begin{lemma}\label{lem classifying space S}\ \\
\vspace{-.15in}
\begin{itemize}
\item[(a)] For any $\alpha\in T(\mu)(\CC)$, $M_\alpha$ is a mixed Hodge structure associated to $\mu$.
\item[(b)] Every mixed Hodge structure associated to $\mu$ is of the form $M_\alpha$ for some $\alpha$ in $T(\mu)(\CC)$.
\item[(c)] For any $\alpha, \beta \in T(\mu)(\CC)$, we have $M_\alpha=M_{\beta}$ if and only if $\alpha^{-1}\beta\in F^0U(\mu)(\CC)$. 
\end{itemize}
\end{lemma}

\begin{proof}
(a) We need to check that the image of $F^\cdot_\alpha W_nM_\CC$ under the quotient map 
\begin{equation}\label{eq33}
W_nM_\CC\ \longrightarrow \ (\widetilde{M}_n)_\CC
\end{equation}
is the Hodge filtration on $(\widetilde{M}_n)_\CC$. We have
\[
F^\cdot_\alpha W_nM_\CC \ = \ \sum\limits_{m\leq n} \alpha_m(F^\cdot (\widetilde{M}_m)_\CC).
\]
Tentatively, let us denote the map Eq. \eqref{eq33} by $\pi_n$. Then the part of the sum above corresponding to $m<n$ goes to zero under $\pi_n$, so that 
\[
\pi_n(F^\cdot_\alpha W_nM_\CC) \ = \ \pi_n \alpha_n(F^\cdot (\widetilde{M}_n)_\CC) \ = \ F^\cdot (\widetilde{M}_n)_\CC.
\]

(b) Let $F^\cdot$ be a filtration on $M_\CC$ such that $M:=(M_\QQ,W_\cdot, F^\cdot)$ is a mixed Hodge structure associated to $\mu$. The map Eq. \eqref{eq33} is the complex underlying map of the natural morphism 
\[W_nM \ \longrightarrow \ \wt{M}_n \, = \, Gr^W_nM,\]
so that by strictness of morphisms of mixed Hodge structures with respect to the Hodge filtration, Eq. \eqref{eq33}
admits a section $\alpha_n$ compatible with the Hodge filtration, i.e. such that
\[
\alpha_n(F^\cdot (\widetilde{M}_n)_\CC) \ \subset \ F^\cdot W_n M_\CC.
\]
Consider $\alpha:=\sum\limits_n  \alpha_n \, \in T(\mu)(\CC)$ (see Remark \ref{elements of T and sections of W_n -> Gr_n}). Then 
\[
F_\alpha^\cdot M_\CC \ = \ \sum\limits_n \alpha_n(F^\cdot (\widetilde{M}_n)_\CC) \ = \ F^\cdot M_\CC.
\]

(c) This is clear from the fact that $F^\cdot_\alpha$ (resp. $F^\cdot_\beta$) was defined by transporting the Hodge filtration of $\wt{M}$ via $\alpha$ (resp. $\beta$) from $\wt{M}_\CC$ to $M_\CC$. 
\end{proof}

The following parametrization of the set of all mixed Hodge structures associated to $\mu$ is immediate from the lemma:

\begin{prop}
The map
\[
T(\mu)(\CC) \ \longrightarrow \ \text{the collection of all mixed Hodge structures associated to $\mu$} 
\]
given by 
\[
\alpha \ \mapsto \ M_\alpha
\]
factors through a bijection
\[
S(\mu)(\CC) \ \longrightarrow \ \text{the collection of all mixed Hodge structures associated to $\mu$}.
\]
\end{prop}

We adopt the following notation for the remainder of the paper. Given any commutative $\CC$-algebra $R$ and any element $\alpha\in T(\mu)(R)$, the image of $\alpha$ in $S(\mu)(R)$ is denoted by $[\alpha]$. Given $s \in S(\mu)(\CC)$, the mixed Hodge structure associated to $\mu$ at $s$ is denoted by $M_s$. The Hodge filtration of $M_s$ is denoted by $F^\cdot_s$ (also denoted by $F^\cdot_\alpha$ if $s=[\alpha]$).

\begin{rem}\label{rem Deligne's splitting}
Recall that given any mixed Hodge structure $N$, Deligne's splitting gives a canonical isomorphism 
\[
a_N \, : \, N_\CC \ \longrightarrow \ Gr^WN_\CC
\]
which preserves the weight and Hodge filtrations (see \cite[Proposition (I.9)]{Mo78} and \cite{De94}). Moreover, $Gr^W a_N$ ( = the map induced by $a_N$ on the weight associated gradeds) is the identity map. Applying this to our setting, for any mixed Hodge structure $N$ associated to the triple $\mu$, we get a canonical element $\alpha_N=a_{N}^{-1}\in T(\mu)(\CC)$ such that $N=M_{\alpha_N}$. In other words, Deligne's splitting gives a set-theoretic section of the map
\[
T(\mu)(\CC) \ \longrightarrow \ S(\mu)(\CC) \, \cong \, \text{the collection of all mixed Hodge structures associated to $\mu$}.
\]
For future use, let us also recall that the isomorphisms $a_N$ above are functorial in $N$, giving an isomorphism of functors 
\[
\omega_0\otimes \CC \ \longrightarrow \ (Gr^W\omega_0)\otimes\CC,
\]
where $\omega_0$ is the forgetful functor from the category of mixed Hodge structures to the category of rational vector spaces.\footnote{In fact, the isomorphism is over $\RR$ (see \cite{De94}), but we will not use this.}
\end{rem}

\subsection{Truncations of the weight filtration of a triple}\label{par: def of fibration}Let $\mu=(M_\QQ,W_\cdot,\wt{M})$ be a triple.\footnote{The material of \S \ref{par: def of fibration} and \S \ref{par: fibers of theta_p} are needed for Theorem \ref{thm1} and the last step of the proof of Theorem \ref{thm2} (\S \ref{sec: properness of the Hodge loci for main thm}).} Fix an integer $p$. The filtration $W_\cdot$ induces filtrations on $W_pM_\QQ$ and $M_\QQ/W_pM_\QQ$; we shall denote these filtrations also by $W_\cdot$. We then have two triples
\[
\mu_p \ := \ (W_pM_\QQ, \, W_\cdot, \, W_p\wt{M})
\]
and 
\[
\mu_{>p} \ := \ (M_\QQ/W_pM_\QQ, \, W_\cdot, \, \wt{M}/W_p\wt{M})
\]
obtained by truncating $\mu$. Every $\alpha\in T(\mu)(R)$ (with $R$ = a commutative $\CC$-algebra) gives elements $\alpha_p\in T(\mu_p)(R)$ and $\alpha_{>p}\in T(\mu_{>p})(R)$ in a natural way (by restriction and passing to the quotient, respectively). We have a morphism
\begin{equation}\label{eq34}
T(\mu) \ \longrightarrow  \ T(\mu_p)\times T(\mu_{>p}) \hspace{.5in} \alpha \, \mapsto \, (\alpha_p, \alpha_{>p}) \, .
\end{equation}
It induces a morphism
\[
\Theta_p \, : \, S(\mu) \ \longrightarrow \ S(\mu_p)\times S(\mu_{>p}) \hspace{.5in} [\alpha] \, \mapsto \, ([\alpha]_p,[\alpha]_{>p}) := ([\alpha_p], [\alpha_{>p}]) \, .
\]
At the level of complex-valued points, thinking of the complex-valued points as parametrizing spaces of mixed Hodge structures associated to the corresponding triples, the map $\Theta_p$ simply sends a mixed Hodge structure $M$ associated to $\mu$ to the pair $(W_pM,M/W_pM)$.

\subsection{Fibers of $\Theta_p$}\label{par: fibers of theta_p}

Let $\mu=(M_\QQ,W_\cdot,\wt{M})$ be a triple and $p\in\ZZ$. The goal of this subsection is to study the fibers of the map $\Theta_p$. 

Let us first fix some notation. Given $x\in S(\mu_p)(\CC)$ (resp. $y\in S(\mu_{>p})(\CC)$), let $(W_pM)_x$ (resp. $(M/W_pM)_y$) be the mixed Hodge structure associated to $\mu_p$ (resp. $\mu_{>p}$) corresponding to $x$ (resp. $y$). These are the mixed Hodge structures with Hodge filtrations $F^\cdot_x$ and $F^\cdot_y$ on $W_pM_\CC$ and $M_\CC/W_pM_\CC$, respectively (the notation for the filtrations being consistent with what was introduced in \S \ref{sec: S parametrizes the MHSs associated to the triple}). 

The set $Sec\, (M_\CC\longrightarrow M_\CC/W_pM_\CC)$ of sections of the map $M_\CC\longrightarrow M_\CC/W_pM_\CC$ is a torsor for $Hom(M_\CC/W_pM_\CC, W_pM_\CC)$ (with the additive action). The algebraification $Sec\, (M_\CC\longrightarrow M_\CC/W_pM_\CC)^a$ (see \S \ref{sec: def of S}) is a torsor for the vector group $Hom(M_\CC/W_pM_\CC, W_pM_\CC)^a$. 

Given
\[
(x,y) \ \in \ S(\mu_p)(\CC) \, \times \, S(\mu_{>p})(\CC),
\]
let 
\[
F^0_{x,y}Hom(M_\CC/W_pM_\CC, W_pM_\CC)
\]
be the subspace of $Hom(M_\CC/W_pM_\CC, W_pM_\CC)$ consisting of linear maps $f$ which are compatible with the filtrations $F^\cdot_y$ and $F^\cdot_x$, by which we mean they satisfy
\[
f(F^\cdot_y(M_\CC/W_pM_\CC)) \ \subset \ F^\cdot_xW_pM_\CC.
\]
The algebraification $F^0_{x,y}Hom(M_\CC/W_pM_\CC, W_pM_\CC)^a$ is the subgroup of $Hom(M_\CC/W_pM_\CC, W_pM_\CC)^a$ whose set of $R$-valued points is the set $F^0_{x,y}Hom((M_\CC/W_pM_\CC)\otimes R, W_pM_\CC\otimes R)$, which is defined to be the set of $R$-linear maps 
\[
(M_\CC/W_pM_\CC)\otimes R \ \longrightarrow \ W_pM_\CC\otimes R
\]
which are compatible with $F^\cdot_y$ and $F^\cdot_x$ (where with abuse of notation, these mean the filtrations induced on the $R$-modules). The quotient
\[
Sec\, (M_\CC\longrightarrow M_\CC/W_pM_\CC)^a/F^0_{x,y}Hom(M_\CC/W_pM_\CC, W_pM_\CC)^a
\] 
is an affine space over $\CC$. Indeed, it is noncanonically isomorphic to 
\[
Hom(M_\CC/W_pM_\CC, W_pM_\CC)^a/F^0_{x,y}Hom(M_\CC/W_pM_\CC, W_pM_\CC)^a.
\] 
Note that 
\begin{equation}\label{eq36}
F^0_{x,y}Hom(M_\CC/W_pM_\CC, W_pM_\CC) \ = \ F^0\inHom((M/W_pM)_y, (W_pM)_x),
\end{equation}
where $\inHom$ is (as everywhere throughout the paper) the internal Hom in the category of mixed Hodge structures and $F^\cdot$ on the right is the Hodge filtration.

Finally, denote the fiber of $\Theta_p$ above $(x,y)$ by $S(\mu)_{(x,y)}$.

\begin{prop}\label{prop: fibers of theta_p}
Let $(x,y)$ be as above.
\begin{itemize}
\item[(a)] There is a canonical isomorphism
\[
Sec \, (M_\CC\longrightarrow M_\CC/W_pM_\CC)^a/F^0_{x,y}Hom(M_\CC/W_pM_\CC, W_pM_\CC)^a \ \longrightarrow \ S(\mu)_{(x,y)}.
\]
\item[(b)] Identifying $S(\mu)(\CC)$ as the set of mixed Hodge structures associated to $\mu$, the image under the isomorphism of Part (a) of a section $\psi$ of $M_\CC\longrightarrow M_\CC/W_pM_\CC$ has Hodge filtration given by
\[
F^\cdot M_\CC \ = \ F^\cdot_x W_pM_\CC \, + \, \psi \, (F^\cdot_y(M_\CC/W_pM_\CC)).
\]
\end{itemize}
\end{prop}

\begin{proof}
We will only include a sketch; the omitted details are all straightforward. First, choose representatives
\[\alpha\in T(\mu_p)(\CC) \hspace{.3in}\text{and}\hspace{.3in} \beta \in T(\mu_{>p})(\CC)\]
for $x$ and $y$, respectively. Given a commutative $\CC$-algebra $R$ and an $R$-linear section $\psi$ of the quotient $M_\CC\otimes R\longrightarrow (M_\CC/W_pM_\CC)\otimes R$, via the canonical isomorphism
\begin{equation}\label{eq23}
\wt{M}_\CC\otimes R \ = \ W_p\wt{M}_\CC\otimes R \ \oplus \ (\wt{M}_\CC/W_p\wt{M}_\CC)\otimes R,
\end{equation}
define 
\[
\gamma_\psi \ := \ \alpha+\psi\beta \, : \, \wt{M}_\CC\otimes R \, \longrightarrow \, M_\CC\otimes R.
\]
Note that, by construction, $\gamma_\psi$ is the unique map $\wt{M}_\CC\otimes R \, \longrightarrow \, M_\CC\otimes R$ fitting in the diagram
\[
\begin{tikzcd}[row sep = large]
  0  \arrow[r] & W_pM_\CC\otimes R   \arrow[r] & M_\CC\otimes R \arrow{r}{\pi} &  (M_\CC/W_pM_\CC)\otimes R \arrow{r} \arrow[bend left=33]{l}[swap]{\psi} & 0 \\
    0 \arrow[r] & W_p\wt{M}_\CC\otimes R \arrow{u}{\alpha}  \arrow[r] & \wt{M}_\CC\otimes R  \arrow{u}{\gamma_\psi} \arrow{r}{} & (\wt{M}_\CC/W_p\wt{M}_\CC)\otimes R  \arrow{u}{\beta} \arrow{r} \arrow[bend left=33]{l}[swap]{\text{canonical}} & 0 \, . 
\end{tikzcd}
\]
Let $T(\mu)_{(\alpha,\beta)}$ be the fiber of Eq. \eqref{eq34} above $(\alpha,\beta)$. We have an isomorphism
\[
Sec \, (M_\CC\longrightarrow M_\CC/W_pM_\CC)^a \ \longrightarrow \ T(\mu)_{(\alpha,\beta)} \hspace{.4in}\psi\mapsto \gamma_\psi \, .
\]
One can check that this isomorphism descends to an isomorphism
\[
Sec \, (M_\CC\longrightarrow M_\CC/W_pM_\CC)^a/F^0_{x,y}Hom(M_\CC/W_pM_\CC, W_pM_\CC)^a \ \longrightarrow \ S(\mu)_{(x,y)}
\]
which does not depend on the choice of $\alpha$ and $\beta$ (this will be the isomorphism of Part (a)). Moreover, note that for a section $\psi$ of $M_\CC\longrightarrow M_\CC/W_pM_\CC$, we have
\begin{align*}
\gamma_\psi(F^\cdot \wt{M}_\CC) \ &= \ \gamma_\psi(F^\cdot W_p\wt{M}_\CC \, + \, F^\cdot (\wt{M}_\CC/W_p\wt{M}_\CC)) \\ 
&= \ \alpha(F^\cdot W_p\wt{M}_\CC) \ + \ \psi\beta F^\cdot (\wt{M}_\CC/W_p\wt{M}_\CC)\\
&= \ F^\cdot_x W_pM_\CC \, + \, \psi \, F^\cdot_y(M_\CC/W_pM_\CC)
\end{align*}
(the first line making use of the canonical isomorphism of Eq. \eqref{eq23}). This shows the formula of Part (b).
\end{proof}

\begin{cor}\label{cor: dim of fibers of theta_p}
Given $x$ and $y$ as above, the fiber $S(\mu)_{(x,y)}$ is an affine space over $\CC$. If $W_pM_\QQ$ and $M_\QQ/W_pM_\QQ$ are both nonzero, then the dimension of $S(\mu)_{(x,y)}$ is at least 1.
\end{cor}

\begin{proof}
The first assertion is immediate from Part (a) of the previous Proposition. The second assertion follows from the same result together with Eq. \eqref{eq36}, on noting that when $W_pM_\QQ$ and $M_\QQ/W_pM_\QQ$ are nonzero, 
\[
\inHom((M/W_pM)_y, (W_pM)_x)
\]
is a nonzero mixed Hodge structure with only negatives weights, and thus its $F^0$ cannot be equal to its entire underlying complex vector space.
\end{proof}

\section{Algebraic families of mixed Hodge structures associated to a triple}\label{sec: algebraic families attached to triples}

\subsection{Definitions and statements}\label{sec: def and statements for alg families}
In this subsection we give the definitions of an algebraic family of mixed Hodge structures associated to a triple and some related notions. We also state two results about these algebraic families which will be needed in \S \ref{sec: proof of main thm} to prove Theorem \ref{thm2}. The proofs of these two results are the subject of the remainder of this section.

Let $\mu=(M_\QQ,W_\cdot,\wt{M})$ be a triple and $X$ a complex variety (reduced scheme of finite type over $\CC$). By an algebraic family of mixed Hodge structures associated to $\mu$ parametrized by $X$ we mean a morphism (of complex varieties)
\[X \ \longrightarrow \ S(\mu).\]

Given an algebraic family as above, for each $x\in X(\CC)$ the image of $x$ corresponds to a mixed Hodge structure associated to $\mu$; we denote this mixed Hodge structure by $M_x$ (with the morphism $X \longrightarrow S(\mu)$ understood from the context). The family 
\[
S(\mu) \ \stackrel{Id}{\longrightarrow} \ S(\mu)
\]
is referred to the universal family associated to $\mu$. Notationally, sometimes instead of writing a family as $X\longrightarrow S(\mu)$, we will write it as $(M_x)_{x\in X(\CC)}$ (with the key being that in that latter notation, $x\mapsto M_x\in S(\mu)(\CC)$ is a morphism).

Given an algebraic family $X \longrightarrow S(\mu)$ of mixed Hodge structures associated to $\mu$ and any subspace
\[
A_\QQ \ \subset \ M_\QQ \, ,
\]
we denote by $X(\CC)_{A_\QQ}\subset X(\CC)$ the set of all $x\in X(\CC)$ such that $A_\QQ$ underlies a mixed Hodge substructure of $M_x$.\footnote{The family with respect to which $X(\CC)_{A_\QQ}$ is defined will be understood from the context, so the fact that the notation $X(\CC)_{A_\QQ}$ does not incorporate $\mu$ or the morphism $X\longrightarrow S(\mu)$ should not lead to confusion.} For any $x\in X(\CC)_{A_\QQ}$, we denote the subobject of $M_x$ supported on $A_\QQ$ by $A_x$. We shall show that:

\begin{prop}\label{lem2 thm2}
Let $X \longrightarrow S(\mu)$ be an algebraic family of mixed Hodge structures associated to $\mu=(M_\QQ,W_\cdot,\wt{M})$. Let $A_\QQ$ be a subspace of $M_\QQ$. Then $X(\CC)_{A_\QQ}$ is a Zariski closed subset of $X(\CC)$.
\end{prop}

With the setting as above, if $X(\CC)_{A_\QQ}$ is nonempty, then there is a subobject of $\wt{M}$ with underlying rational vector space
\[
\wt{A}_\QQ \ := \ Gr^WA_\QQ \ \subset \ \wt{M}_\QQ \, . 
\]
We shall denote this subobject by $\wt{A}$. Note that for every $x\in X(\CC)_{A_\QQ}$ we then have
\[
Gr^WA_x \ = \ \wt{A},
\]
and the mixed Hodge structures $A_x$ and $M_x/A_x$ are associated to the triples 
\[
(A_\QQ,W_\cdot,\widetilde{A})
\]
and
\[(M_\QQ/A_\QQ, W_\cdot,\widetilde{M}/\widetilde{A})\]
respectively (where with abuse of notation, $W_\cdot$ denotes the induced weight filtrations on $A_\QQ$ and $M_\QQ/A_\QQ$).

Let us say $A_\QQ$ is a global Hodge subspace with respect to the family if 
\[
X(\CC)_{A_\QQ} \ = \ X(\CC),
\]
i.e. if $A_\QQ$ underlies a subobject of $M_x$ for every $x\in X(\CC)$. Our next result is the following:

\begin{prop}\label{prop quotients for alg families}
Let $X \longrightarrow S(\mu)$ be an algebraic family of mixed Hodge structures associated to $\mu=(M_\QQ,W_\cdot,\wt{M})$. Let $A_\QQ\subset M_\QQ$ be a global Hodge subspace with respect to the family. Then we have algebraic families $(A_x)_{x\in X(\CC)}$ and $(N_x/A_x)_{x\in X(\CC)}$, parametrized by $X$, respectively associated to the triples $\mu |_{A_\QQ}:=(A_\QQ,W_\cdot,\widetilde{A})$ and $\mu/A_\QQ:=(M_\QQ/A_\QQ, W_\cdot,\widetilde{M}/\widetilde{A})$.
\end{prop}

The remainder of this section is devoted to the proofs of the two propositions above.

\subsection{A useful lemma}
Given a rational vector space $M_\QQ$ and a complex subspace $A_\CC$ of $M_\CC:=M_\QQ\otimes\CC$, by saying $A_\CC$ is defined over $\QQ$ we mean that 
\[
A_\CC \ = \ A_\QQ\otimes \CC
\]
for some subspace $A_\QQ\subset M_\QQ$. By basic linear algebra, (i) if $A_\QQ$ exists then it is unique and is equal to $A_\CC\cap M_\QQ$, and (ii) $A_\CC$ is defined over $\QQ$ if and only if 
\[
A_\CC \ = \ (A_\CC\cap M_\QQ)\otimes \CC.
\]

The following lemma will be used in proving Proposition \ref{lem2 thm2}; the lemma tells us when a subobject of the associated graded of a mixed Hodge structure $M$ lifts to a subobject of $M$:

\begin{lemma}\label{lem1 thm2} Let $M$ be a mixed Hodge structure, and $\widetilde{A}$ a Hodge substructure of $\wt{M}:=Gr^WM$. Then there exists a subobject $A$ of $M$ such that $Gr^WA=\wt{A}$ if and only if there exists an isomorphism
\[
\alpha \; : \; \wt{M}_\CC \ \longrightarrow \ M_\CC
\]
which satisfies the following properties:
\begin{itemize}
\item[(i)] $\alpha$ is compatible with the weight filtrations and $Gr^W\alpha $ is the identity map on $\wt{M}_\CC$.
\item[(ii)] $\alpha$ is compatible with the Hodge filtrations.
\item[(iii)] For every $n$, the subspace $\alpha(W_n\wt{A}_\CC)$ of $W_nM_\CC$ is defined over $\QQ$.
\end{itemize}
Moreover, if $A$ exists, then it is unique, and for any $\alpha$ as above, we have $A_\CC=\alpha(\wt{A}_\CC)$.
\end{lemma}

\begin{proof}
We first discuss uniqueness. One way to see this is as follows: Let $\omega_0$ be the forgetful functor from $\mathbf{MHS}$ to rational vector spaces. Then $\omega^{Gr}:=\omega_0\circ Gr^W$ is also a fiber functor (sending $V\mapsto Gr^WV_\QQ$ for any mixed Hodge structure $V$). Denote (tentatively) the Tannakian fundamental group of an object $V$ with respect to $\omega^{Gr}$ by $\mathcal{G}(V)$. If there were objects $A_1$ and $A_2$ with associated graded $\wt{A}$,  then
\[\omega^{Gr}A_1=\omega^{Gr}A_2=\wt{A}_\QQ.\]
The representation of $\mathcal{G}(M)$ on $\wt{A}_\QQ$ corresponding to each $A_i$ is the sub-representation of the action of $\mathcal{G}(M)$ on $\wt{M}_\QQ$ corresponding to $M$, so that the two actions on $\wt{A}_\QQ$ corresponding to $A_1$ and $A_2$ are the same. The identity map on $\wt{A}_\QQ$ gives a $\mathcal{G}(M)$-equivariant isomorphism 
\[
\omega^{Gr}A_1 \ \longrightarrow \ \omega^{Gr}A_2,
\] 
where the action on $\omega^{Gr}A_i$ corresponds to $A_i$. Since $\omega^{Gr}$ is an equivalence of categories from the Tannakian category generated by $M$ to the category of finite-dimensional representations of $\mathcal{G}(M)$, this gives an isomorphism
\[
A_1 \ \longrightarrow \ A_2
\] 
which is compatible with the inclusions into $M$. Thus $A_1 =A_2$. 

\underline{``If" implication}: Given $\alpha$ as in the statement, set
\begin{align*}
A_\CC \ &:= \ \alpha(\wt{A}_\CC)\\
A_\QQ \ &:= \ A_\CC\cap M_\QQ.
\end{align*}
For $R\in \{\QQ,\CC\}$, define the weight filtration on $A_R$ by 
\[W_\cdot A_R \ := \ A_R\cap W_\cdot M_R\]
By (i) we have $W_\cdot A_\CC=\alpha(W_\cdot \wt{A}_\CC)$. By (iii), 
\[
W_\cdot A_\CC \ = \ (W_\cdot A_\QQ)\otimes \CC.
\]
The restriction
\[
\alpha |_{\wt{A}_\CC} \, : \, \wt{A}_\CC \ \longrightarrow \ A_\CC,
\]
of $\alpha$ to $\wt{A}_\CC$ is compatible with weight filtrations. Since (being a subobject of a graded Hodge structure) $\wt{A}$ is graded by weight, we have $\wt{A}\cong Gr^W\wt{A}$. Applying $Gr^W$ to
\begin{equation}\label{eq4}
\begin{tikzcd}[row sep = tiny]
  \wt{A}_\CC  \arrow{r}{ \alpha |_{\wt{A}_\CC}}[swap]{\simeq} & A_\CC \\
  \bigcap & \bigcap \\
    \wt{M}_\CC \arrow{r}{ \ \ \alpha \ \ }[swap]{\simeq} & M_\CC 
\end{tikzcd}
\end{equation}
and on recalling that $Gr^W\alpha$ is the identity, we see that $\wt{A}_\CC=Gr^WA_\CC$, and that after passing to associated gradeds $\alpha |_{\wt{A}_\CC}$ is the identity map. Use $\alpha |_{\wt{A}_\CC}$ to transport the Hodge filtration on $\wt{A}_\CC$ to obtain a filtration $F^\cdot _{\alpha |_{\wt{A}_\CC}}$ (notation consistent with \S \ref{sec: S parametrizes the MHSs associated to the triple}) on $A_\CC$. Then by Lemma \ref{lem classifying space S}(a),
\[
A \ := \ (A_\QQ, W_\cdot , F^\cdot _{\alpha |_{\wt{A}_\CC}})
\]
is a mixed Hodge structure with $Gr^WA=\wt{A}$. The inclusion $\iota: A_\QQ\subset M_\QQ$ is a morphism of mixed Hodge structures. (Compatibility with the weight filtration is by definition. As for the Hodge filtration, the left inclusion in Eq. \eqref{eq4} is compatible with the Hodge filtrations, and the Hodge filtrations on $A_\CC$ and $M_\CC$ are transports of the Hodge filtrations on $\wt{A}_\CC$ and $\wt{N}_\CC$ via the horizontal isomorphisms.)

\underline{``Only if" implication}: This follows from functoriality of Deligne's splitting. The isomorphism 
\[\alpha=a_M^{-1}  \, : \,  \wt{M}_\CC \ \longrightarrow \ M_\CC\]
given by Deligne's splitting does the job (see Remark \ref{rem Deligne's splitting}).

The very last assertion in the statement follows from uniqueness and the proof of the ``if" implication.

\end{proof}

With the setting as in the previous lemma, if $A$ exists we say $\wt{A}$ lifts to $A$. 

\subsection{Proof of Proposition \ref{lem2 thm2}}\label{sec: proof of closed Hodge loci result} We now prove Proposition \ref{lem2 thm2}. It suffices to prove the statement for the universal family, that is, we need to show that the subset $S(\mu)(\CC)_{A_\QQ}$ of $S(\mu)(\CC)$ consisting of all $s\in S(\mu)(\CC)$ such that $A_\QQ$ underlies a subobject of $M_s$ is Zariski closed. We may assume that $S(\mu)(\CC)_{A_\QQ}$ is nonempty. Note that  
\[
\wt{A} \ := \ Gr^W A_s \ \subset \ \wt{M}
\]
is independent of $s\in S(\mu)(\CC)_{A_\QQ}$, as it is the subobject of $\wt{M}$ with underlying rational vector space $Gr^WA_\QQ\subset \wt{M}_\QQ$. The set $S(\mu)(\CC)_{A_\QQ}$ is the set of all $s$ in $S(\mu)(\CC)$ such that $\wt{A}$ lifts to a subobject of $M_s$, and moreover the underlying complex vector space of this subobject is $A_\CC :=A_\QQ\otimes \CC$.

Since $S(\mu)(\CC)_{A_\QQ}$ is nonempty, for every $n$, the vector space $W_nA_\CC:= A_\CC\cap W_n M_\CC$ is defined over $\QQ$. By Lemma \ref{lem1 thm2}, a point $s$ in $S(\mu)(\CC)$ belongs to $S(\mu)(\CC)_{A_\QQ}$ if and only if there exists $\alpha\in T(\mu)(\CC)$ with $s=[\alpha]$ such that $\alpha(\wt{A}_\CC)=A_\CC$. 

Let $H$ be the subgroup of $GL(M_\CC)$ whose group of $\CC$-valued points consists of the elements of $GL(M_\CC)(\CC)$ which preserve the weight filtration, induce identity on $Gr^WM_\CC$, and stabilize $A_\CC$. Then $H$ acts on $T(\mu)$ on the left by post-composition. This action induces an action of $H$ on
\[
S(\mu) \ = \ T(\mu)/F^0U(\mu).
\]
By the characterization of the elements of $S(\mu)(\CC)_{A_\QQ}$ described above, this set is an orbit of the action of $H(\CC)$ on $S(\mu)(\CC)$. By the theorem of Kostant-Rosenlicht (see for instance, \cite[Theorem 17.64]{Mi17}), such orbits are all Zariski closed. (Note that $H$ is unipotent and $S(\mu)$ is affine.)

\begin{rem}
Proposition \ref{lem2 thm2} together with uniqueness of the lifting in Lemma \ref{lem1 thm2} has the following direct consequence. Let $X\longrightarrow S(\mu)$ be an algebraic family of mixed Hodge structures associated to $\mu=(M_\QQ,W_\cdot, \wt{M})$ and $\wt{A}$ a subobject of $\wt{M}$. Let $X(\CC)^{\wt{A}}$ be the set of all $x\in X(\CC)$ such that $\wt{A}$ lifts to a subobject of $M_x$. For each $x\in X(\CC)^{\wt{A}}$ denote the lift by $A_x$. Then one of the following statements is true:
\begin{itemize}
\item[(i)] $X(\CC)^{\wt{A}}=X(\CC)$, with the underlying rational vector space of $A_x$ constant throughout the family.
\item[(ii)] $X(\CC)^{\wt{A}}$ is the union of countably many disjoint proper Zariski closed subsets of $X(\CC)$.
\end{itemize}
Indeed, $X(\CC)^{\wt{A}}=\bigsqcup\limits_{A_\QQ} X(\CC)_{A_\QQ}$ with the union taken over all subspaces $A_\QQ$ of $M_\QQ$ such that $Gr^WA_\QQ=\wt{A}_\QQ$. 
\end{rem}

\subsection{Proof of Proposition \ref{prop quotients for alg families}}
The proof of Proposition \ref{prop quotients for alg families} relies on the following fact about torsors for unipotent algebraic groups in characteristic zero: if $K$ is a unipotent group over a field $k$ of characteristic zero, $V$ is an affine variety over $k$, then every $K$-torsor over $V$ is trivial, in the sense that it admits a section (see \cite[Proposition 16.55]{Mi17}). 

We now give the proof of Proposition \ref{prop quotients for alg families}. It is enough to prove the assertion for the family
\[
S(\mu)_{A_\QQ} \ \hookrightarrow  \ S(\mu),
\]
where $S(\mu)_{A_\QQ}$ as the subvariety of $S(\mu)$ whose set of closed points is $S(\mu)(\CC)_{A_\QQ}$, the latter set being defined with respect to the universal family, as in the proof of Proposition \ref{lem2 thm2} (\S \ref{sec: proof of closed Hodge loci result}). Let $H$ also be as in the proof of Proposition \ref{lem2 thm2}. Then $S(\mu)_{A_\QQ}$ is an orbit of the action of $H$ on $S(\mu)$. We are assuming that $S(\mu)_{A_\QQ}$ is nonempty. Note that $S(\mu)_{A_\QQ}$ is affine because it is a closed subvariety of $S(\mu)$. 

Fix $s_0\in S(\mu)(\CC)_{A_\QQ}$. By Lemma \ref{lem1 thm2}, we can write $s_0$ as $[\alpha_0]$ for some $\alpha_0\in T(\mu)(\CC)$ with $\alpha_0(\wt{A}_\CC)=A_\CC$. Let $K$ be the stabilizer of $s_0$ in $H$ for the action of $H$ on $S(\mu)$. The orbit morphism
\[
\pi \, : \, H \ \longrightarrow S(\mu)_{A_\QQ} \hspace{.5in} h\mapsto h\cdot s_0
\]
is a $K$-torsor over the affine variety $S(\mu)_{A_\QQ}$. Thus $\pi$ admits a section $\beta$.

We have a morphism
\[
(- \circ \alpha_0)\, |_{\wt{A}_\CC} \, : \, H \ \longrightarrow \ T(\mu |_{A_\QQ}) \hspace{.5in} h \, \mapsto \, h  \alpha_0\, |_{\wt{A}_\CC} \, ,
\] 
where given $h\in H(R)$, by $h  \alpha_0 \, |_{\wt{A}_\CC}$ we mean the restriction of $h  \alpha_0$ to an isomorphism $\wt{A}_\CC\otimes R\longrightarrow A_\CC\otimes R$. We also have a morphism 
\[
\overline{-\circ \alpha_0} \, : \, H \ \longrightarrow \ T(\mu/A_\QQ) \hspace{.5in} h \, \mapsto \, \overline{h  \alpha_0} \, ,
\]
where for every $h\in H(R)$, by $\overline{h  \alpha_0}$ we denote the isomorphism $(\wt{M}_\CC/\wt{A}_\CC)\otimes R\longrightarrow (M_\CC/A_\CC)\otimes R$ induced by $h  \alpha_0$. The compositions
\[
S(\mu)_{A_\QQ} \ \stackrel{\beta}{\longrightarrow} \ H \ \xrightarrow{(- \circ \alpha_0)\, |_{\wt{A}_\CC}} \ T(\mu |_{A_\QQ}) \ \longrightarrow \ S(\mu |_{A_\QQ})
\]
and
\[
S(\mu)_{A_\QQ} \ \stackrel{\beta}{\longrightarrow} \ H \ \xrightarrow{\overline{-\circ \alpha_0}} \ T(\mu/A_\QQ) \ \longrightarrow \ S(\mu/A_\QQ)
\]
are the desired algebraic families, as is easily seen by looking at the complex points. (In particular, the compositions are independent of the choice of $\beta$.)


\section{Proof of Theorem \ref{thm2}}\label{sec: proof of main thm}

\subsection{Setup and reduction to the case of global Hodge subspaces}\label{par: setup for proof of main thm}
The goal of this section is to use the results of the previous section and Corollary 5.3.2 of \cite{EM21b} to deduce Theorem \ref{thm2}. Throughout the section, we shall fix the integer $p$ and the triple $\mu=(M_\QQ,W_\cdot, \wt{M})$. We assume that both $W_pM_\QQ$ and $M_\QQ/W_pM_\QQ$ are nonzero, as otherwise $D^p(\mu)$ is empty (see the Introduction for the definition of $D^p(\mu)$). For any proper subspace $A_\QQ\subset Hom(M_\QQ/W_pM_\QQ, W_pM_\QQ)$, let $D^{p}(\mu)_{A_\QQ}$ be the subset of $S(\mu)(\CC)$ consisting of all $s\in S(\mu)(\CC)$ such that 
the following two conditions hold:
\begin{itemize}
\item[(i)] $A_\QQ$ underlies a subobject of $\inHom(M_s/W_pM_s, W_pM_s)$.
\item[(ii)] $\fu_p(M_s)\subset A_\QQ$.
\end{itemize}
We then have
\[
D^{p}(\mu) \ = \  \bigcup\limits_{A_\QQ}\,  D^{p}(\mu)_{A_\QQ},
\]
where the union is over all proper subspaces $A_\QQ$ of $Hom(M_\QQ/W_pM_\QQ, W_pM_\QQ)$. Thus it suffices to show that each $D^{p}(\mu)_{A_\QQ}$ is contained in a union of countably many proper Zariski closed subsets of $S(\mu)(\CC)$.

We will be working with the following two triples:
\begin{align*}
\nu  \ &:= \ \bigm(Hom(M_\QQ/W_pM_\QQ, W_pM_\QQ) , \, W_\cdot, \, \inHom(\wt{M}/W_p\wt{M}, W_p\wt{M})\bigm)\\
\nu' \ &:= \ \bigm(Hom(M_\QQ/W_pM_\QQ, M_\QQ)^\dagger, \, W_\cdot, \, \inHom(\wt{M}/W_p\wt{M}, \wt{M})^\dagger\bigm)
\end{align*}
We refer the reader to \S \ref{sec: method of proofs in Intro} for the definition of $\inHom(N/W_pN, N)^\dagger$ and its underlying rational vector space $Hom(N_\QQ/W_pN_\QQ, N_\QQ)^\dagger$ for any mixed Hodge structure $N$. The weight filtrations in the triples $\nu$ and $\nu'$ are those induced by the weight filtration on $M_\QQ$. (For $\nu'$, consider the weight filtration on $Hom(M_\QQ/W_pM_\QQ, M_\QQ)$ first and then restrict it to the daggered subspace.)

\begin{lemma}\label{alg families coming from linear algebra}
We have algebraic families
\[
(\inHom(M_s/W_pM_s, W_pM_s))_{s\in S(\mu)(\CC)}
\]
and 
\[
(\inHom(M_s/W_pM_s, M_s)^\dagger)_{s\in S(\mu)(\CC)}
\]
of mixed Hodge structures associated respectively to the triples $\nu$ and $\nu'$. 
\end{lemma}
\begin{proof}
For a filtered module $(V,W_\cdot)$ over a commutative algebra $k$ over $\QQ$, to simplify the notation let us use the canonical injections to consider both $Hom(V/W_pV,W_pV)$ and $Hom(V/W_pV,V)^\dagger$ (the latter defined in a way analogous to the case of $k=\QQ$) as subspaces of $End(V)$. The desired morphisms $S(\mu)\longrightarrow S(\nu)$ and $S(\mu)\longrightarrow S(\nu')$ descend from morphisms $T(\mu)\longrightarrow T(\nu)$ and $T(\mu)\longrightarrow T(\nu')$ given by conjugation, i.e. sending $\alpha\in T(\mu)(R)$ (for a commutative $\CC$-algebra $R$) to conjugation by $\alpha$.
\end{proof}

Let $A_\QQ\subsetneq Hom(M_\QQ/W_pM_\QQ, W_pM_\QQ)$. Applying Proposition \ref{lem2 thm2} to the families  above, the set
\[
S(\mu)(\CC)_{A_\QQ} \ = \ \{s\in S(\mu)(\CC) \, : \,  \text{$A_\QQ$ underlies a subobject of $\inHom(M_s/W_pM_s, W_pM_s)$}\}\footnote{While in the definition of $S(\mu)(\CC)_{A_\QQ}$ we worked with the family $(\inHom(M_s/W_pM_s, W_pM_s))_{s\in S(\mu)(\CC)}$, it is clear that working with the family $(\inHom(M_s/W_pM_s, M_s)^\dagger)_{s\in S(\mu)(\CC)}$ would result in the same $S(\mu)(\CC)_{A_\QQ}$.}
\]
is a Zariski closed subset of $S(\mu)(\CC)$. Since $D^{p}(\mu)_{A_\QQ}$ is contained in $S(\mu)(\CC)_{A_\QQ}$, if $S(\mu)(\CC)_{A_\QQ} \neq S(\mu)(\CC)$ we are done. It remains to consider the case where $S(\mu)_{A_\QQ} =S(\mu)(\CC)$.

\subsection{Application of a result of \cite{EM21b} - Reduction to the study of Hodge loci in quotient families}
Assume that $A_\QQ$ is a proper subspace of $Hom(M_\QQ/W_pM_\QQ, W_pM_\QQ)$ with $S(\mu)(\CC)_{A_\QQ} =S(\mu)(\CC)$. In what follows by a Hodge class in a rational mixed Hodge structure $N$ we mean a Hodge class of weight zero, i.e. 
an element of
\[
W_0N_\QQ \, \cap \, F^0N_\CC. 
\]
The data of a Hodge class $w$ in $N$ is equivalent to the data of a morphism $\mathbbm{1}\longrightarrow N$ (passage from the former to the latter: define the morphism by sending $1\mapsto w$; passage from the latter to the former: take $w$ to be the image of $1\in\QQ$). 

Fix $s\in S(\mu)(\CC)$ for the time being. With notation as in \S \ref{sec: def and statements for alg families}, $A_s$ denotes the subobject of $\inHom(M_s/W_pM_s, W_pM_s)$ (and $\inHom(M_s/W_pM_s, M_s)^\dagger$) supported on $A_\QQ$. Consider the pushforward $\mathcal{E}_p(M_s)/A_s$ of the extension 
\[
\mathcal{E}_p(M_s): \quad 0 \ \longrightarrow \ \inHom(M_s/W_pM_s, W_pM_s) \longrightarrow \ \inHom(M_s/W_pM_s, M_s)^\dagger \stackrel{\lambda}{\longrightarrow} \ \mathbbm{1} \longrightarrow \ 0
\]
(see \S \ref{sec: method of proofs in Intro}) along the quotient map
\[
\inHom(M_s/W_pM_s, W_pM_s) \ \longrightarrow \ \inHom(M_s/W_pM_s, W_pM_s)/A_s \, .
\]
Under Hypotheses\footnote{This is the only place in the argument where these hypotheses come into play.} (i) and (ii) of Theorem \ref{thm2}, Corollary 5.3.2 of \cite{EM21b} tells us that $\mathcal{E}_p(M_s)/A_s$ splits if and only if $\fu_p(M_s)$ is contained in $A_\QQ$ (or in other words, $\underline{\fu}_p(M_s)\subset A_s$, see \S \ref{sec: method of proofs in Intro} for definition of $\underline{\fu}_p(M_s)$). The extension $\mathcal{E}_p(M_s)/A_s$ is given by
\[
0 \ \longrightarrow \ \inHom(M_s/W_pM_s, W_pM_s)/A_s \longrightarrow \ \inHom(M_s/W_pM_s, M_s)^\dagger/A_s \stackrel{\overline{\lambda}}{\longrightarrow} \ \mathbbm{1} \longrightarrow \ 0
\]
($\overline{\lambda}$ being induced by $\lambda$). There is a bijection 
\begin{align*}
\text{sections $\gamma$ of $\overline{\lambda}$ in the category} \ \ \ \ & \longleftrightarrow \ \ \ \ \text{Hodge classes $w$ in $\inHom(M_s/W_pM_s, M_s)^\dagger/A_s$} \\
\text{of mixed Hodge structures} \ \ \ \ \ \ \ & \ \ \ \ \ \ \ \ \  \ \ \text{with $\overline{\lambda}(w)=1$}\\
\gamma \ \ \ & \ \leftrightarrow \ \ \ \ \ \gamma(1) \, .
\end{align*}
Note that since the weights of $\inHom(M_s/W_pM_s, W_pM_s)/A_s$ are all negative, $\overline{\lambda}$ admits at most one section in the category of mixed Hodge structures. 

For any element $w$ of
\[
Hom(M_\QQ/W_pM_\QQ, M_\QQ)^\dagger/A_\QQ 
\]
with\footnote{Note that as maps of vector spaces the maps $\lambda$ and $\overline{\lambda}$ do not depend on $s$.} $\overline{\lambda}(w)=1$, consider the Hodge locus
\[
D^{p}(\mu)_{A_\QQ,w} \ := \ \{s\in S(\mu)(\CC) \, : \, \text{$w$ is a Hodge class in $\inHom(M_s/W_pM_s, M_s)^\dagger/A_s$}\}.
\]
Then combining the above observations, we have 
\[
D^{p}(\mu)_{A_\QQ} \ = \ \bigsqcup\limits_{w} \, D^{p}(\mu)_{A_\QQ,w}
\]
where $w$ runs through all $w$ as above. We will be done if we show that each $D^p(\mu)_{A_\QQ,w}$ is a proper Zariski closed subset of $S(\mu)(\CC)$. 

\subsection{Zariski closedness of each Hodge locus} 
Zariski closedness of each $D^{p}(\mu)_{A_\QQ,w}$ (with $A_\QQ$ and $w$ as above) follows from Propositions \ref{lem2 thm2} and \ref{prop quotients for alg families}. Indeed, $D^{p}(\mu)_{A_\QQ,w}$ is, in the notation of \S \ref{sec: def and statements for alg families}, just equal to $S(\mu)(\CC)_{span\{w\}}$ for the algebraic family 
\[(\inHom(M_s/W_pM_s, M_s)^\dagger/A_s)_{s\in S(\mu)(\CC)}\, .\] 
(Note that if $\text{span}\{w\}$ is a subobject of $\inHom(M_s/W_pM_s, M_s)^\dagger/A_s$, then this subobject will have to be a copy of $\mathbbm{1}$: being a 1-dimensional subobject of an object with weights $\leq 0$, it will be a copy of $\QQ(n)$ for some $n\geq 0$, and since $\overline{\lambda}(w)\neq 0$ we must have $n=0$.)

\subsection{Properness of each Hodge locus}\label{sec: properness of the Hodge loci for main thm}
The final step is to show that each Hodge locus $D^{p}(\mu)_{A_\QQ,w}$ is a proper subset of $S(\mu)(\CC)$. In other words, given any $f\in Hom(M_\QQ/W_pM_\QQ, M_\QQ)^\dagger$ with $\lambda(f)=1$ (i.e. given any section of $M_\QQ\longrightarrow M_\QQ/W_pM_\QQ$), the goal is to show that there is $s\in S(\mu)(\CC)$ such that 
\[f+A_\QQ \ \in \ Hom(M_\QQ/W_pM_\QQ, M_\QQ)^\dagger/A_\QQ\]
is not a Hodge class in $\inHom(M_s/W_pM_s, M_s)^\dagger/A_s$. 

Recall from \S \ref{par: def of fibration} and \S \ref{par: fibers of theta_p} that we have a map
\[
\Theta_p \, : \, S(\mu)(\CC) \ \longrightarrow \ S(\mu_p)(\CC) \times S(\mu_{>p})(\CC) \hspace{.3in} s\mapsto (s_p, s_{>p})
\]
(here considered at the level of complex points). For any 
\[
(x,y) \ \in \ S(\mu_p)(\CC) \times S(\mu_{>p})(\CC),
\]
the fiber $S(\mu)(\CC)_{(x,y)}$ of $\Theta_p$ above $(x,y)$ is canonically isomorphic to the set of sections of the linear map $M_\CC\longrightarrow M_\CC/W_pM_\CC$, modded out by the additive action of 
\[F^0\inHom((M/W_pM)_y, (W_pM)_x)\]
(see Proposition \ref{prop: fibers of theta_p}). Denote the element of $S(\mu)(\CC)_{(x,y)}$ corresponding to a section $\psi$ by $[\psi]$. The Hodge filtration 
of $M_{[\psi]}$ is given by the formula
\begin{equation}\label{eq37}
F^\cdot_{[\psi]} M_\CC \ = \ F^\cdot_xW_pM_\CC \, + \, \psi(F^\cdot_y (M_\CC/W_pM_\CC)).
\end{equation}

We shall show that for any $f\in Hom(M_\QQ/W_pM_\QQ, M_\QQ)^\dagger$ with $\lambda(f)=1$ and any $(x,y)$ in $S(\mu_p)(\CC) \times S(\mu_{>p})(\CC)$, there is a point $s$ in the fiber $S(\mu)(\CC)_{(x,y)}$ such that $f+A_\QQ$ is not a Hodge class in $\inHom(M_s/W_pM_s, M_s)^\dagger/A_s$. More precisely, we will see that:
\begin{prop}\label{prop: properness of each loci in thm 2}
Let $f\in Hom(M_\QQ/W_pM_\QQ, M_\QQ)^\dagger$ with $\lambda(f)=1$. Fix arbitrary $(x,y) \in S(\mu_p)(\CC) \times S(\mu_{>p})(\CC)$. Let $A_\QQ$ be a proper subspace of $Hom(M_\QQ/W_pM_\QQ, W_pM_\QQ)$ which underlines a mixed Hodge substructure of $\inHom((M/W_pM)_y, (W_pM)_x)$.\footnote{This is guaranteed if $S(\mu)(\CC)_{A_\QQ} =S(\mu)(\CC)$. Indeed, since $W_pM_s=(W_pM)_{s_p}$ and $M_s/W_pM_s=(M/W_pM)_{s_{>p}}$, the subspace $A_\QQ$ underlines a subobject of $\inHom((M/W_pM)_y, (W_pM)_x)$ if and only if it underlines a subobject of $\inHom(M_s/W_pM_s, W_pM_s)$ for one (and hence every) $s\in  S(\mu)(\CC)_{(x,y)}$.} Then there is $s\in S(\mu)(\CC)_{(x,y)}$ such that $f+A_\QQ$ is not a Hodge class in $\inHom(M_s/W_pM_s, M_s)^\dagger/A_s$. 
\end{prop}

Let us first give two lemmas:

\begin{lemma}\label{lem: part 2 properness (easy but important observation)}
Fix $(x,y) \in S(\mu_p)(\CC) \times S(\mu_{>p})(\CC)$. Let $f$ and $\psi$ be sections of $M_\CC\longrightarrow M_\CC/W_pM_\CC$. Let $s=[\psi] \in S(\mu)(\CC)_{(x,y)}$. Then $f \in F^0\inHom(M_s/W_pM_s, M_s)^\dagger$ if and only if 
\[\psi \ \in \ f+ F^0\inHom((M/W_pM)_y, (W_pM)_x).\]
\end{lemma}

\begin{proof}
In view of the definition of the Hodge filtration of an internal Hom object, to say $f \in F^0\inHom(M_s/W_pM_s, M_s)^\dagger$ is to say that $f$ satisfies $f(F^\cdot_{y}(M_\CC/W_pM_\CC)) \subset F^\cdot_{s}M_\CC$. The latter statement is by Eq. \eqref{eq37} (applied to the section $f$) equivalent to $F^\cdot_{[f]}=F^\cdot_s$, or in other words, to $[f]=s$. The claim follows.
\end{proof}

The next lemma builds on the previous one:
\begin{lemma}\label{lem 2 properness of Hodge loci}
Let $f$, $(x,y)$ and $A_\QQ$ be as in Proposition \ref{prop: properness of each loci in thm 2}. Let $s=[\psi]\in S(\mu)(\CC)_{(x,y)}$, with $\psi$ a section of $M_\CC\longrightarrow M_\CC/W_pM_\CC$. Then $f+A_\QQ$ is a Hodge class in the quotient $\inHom(M_s/W_pM_s, M_s)^\dagger/A_s$ if and only if 
\[\psi \ \in \ f+A_\CC+F^0\inHom((M/W_pM)_y, (W_pM)_x).\]
\end{lemma}

\begin{proof}
The element $f+A_\QQ$ is a Hodge class in $\inHom(M_s/W_pM_s, M_s)^\dagger/A_s$ if and only if it is in $F^0$ of the latter mixed Hodge structure. By strictness of morphisms of mixed Hodge structures with respect to the Hodge filtration, this is equivalent to 
\[(f+A_\CC) \, \cap \, F^0\inHom(M_s/W_pM_s, M_s)^\dagger\]
being nonempty. Combining with Lemma \ref{lem: part 2 properness (easy but important observation)}, we get the claim.
\end{proof}

We are ready to deduce Proposition \ref{prop: properness of each loci in thm 2}, and hence complete the proof of Theorem \ref{thm2}. 

\begin{proof}[Proof of Proposition \ref{prop: properness of each loci in thm 2}]
With notation as in the statement of the Proposition, suppose that $f+A_\QQ$ is a Hodge class in $\inHom(M_s/W_pM_s, M_s)^\dagger/A_s$ for every $s\in S(\mu)(\CC)_{(x,y)}$. By Lemma \ref{lem 2 properness of Hodge loci}, 
\[ f+A_\CC+F^0\inHom((M/W_pM)_y, (W_pM)_x)\] 
contains (and hence is equal to) the set of all sections of $M_\CC\longrightarrow M_\CC/W_pM_\CC$. Thus 
\[A_\CC+F^0\inHom((M/W_pM)_y,(W_pM)_x) \, = \, Hom(M_\CC/W_pM_\CC, W_pM_\CC).\]
It follows that
\[F^0\, \left(\inHom((M/W_pM)_y,(W_pM)_x)/ A \right) \ = \ \inHom((M/W_pM)_y,(W_pM)_x)/A \, ,\]
where $A$ is the subobject of $\inHom((M/W_pM)_y,(W_pM)_x)$ supported on $A_\QQ$. But this is absurd since $\inHom((M/W_pM)_y,(W_pM)_x)/A$ is a nonzero mixed Hodge structure of only negative weights. 
\end{proof}

\section{Proof of Theorem \ref{thm1}}\label{sec: proof of Thm C}

The goal of this section is to prove Theorem \ref{thm1}. Throughout, we fix the triple $\mu$, an integer $p$, and a point $(x,y)\in S(\mu_p)(\CC)\times S(\mu_{>p})(\CC)$. The theorem is trivial if $M_\CC/W_pM_\CC$ or $W_pM_\CC$ is zero, so again we will assume that these are both nonzero. Notations such as $(W_pM)_x, (M/W_pM)_y$, $\lambda$, etc. mean the same as before. Thus for every $s\in S(\mu)(\CC)_{(x,y)}$, we have $W_pM_s=(W_pM)_x$ and $M_s/W_pM_s=(M/W_pM)_y$, so that 
\[
\underline{\fu}_p(M_s) \ \subset \ \inHom((M/W_pM)_y, (W_pM)_x).
\]
(See \S \ref{sec: method of proofs in Intro} for the definition of $\underline{\fu}_p(M_s)$.)

\subsection{Setup}\label{par: setup proof of Theorem C} For every proper subobject 
\[ A \ \subset \ \inHom((M/W_pM)_y, (W_pM)_x), \]
set
\[
{D'}\, ^p(\mu)_{A}  \ := \ \{s\in S(\mu)(\CC)_{(x,y)} \, : \,  \underline{\fu}_p(M_s)\subset A\}.\footnote{In the notation of \S \ref{sec: proof of main thm}, this is the intersection of $D^p(\mu)_{A_\QQ}$ and $S(\mu)(\CC)_{(x,y)}$. We used $A$ instead of $A_\QQ$ in the notation because $A$ is a fixed mixed Hodge structure when $x$ and $y$ are fixed.} 
\]
Since $\inHom((M/W_pM)_y, (W_pM)_x)$ has only countably many subobjects, it is enough to show that for each of those proper subobjects $A$, the set $D'\, ^p(\mu)_{A}$ above is a union of countably many proper Zariski closed subsets of $S(\mu)(\CC)_{(x,y)}$. By Theorem 3.3.1 of \cite{EM21a} (which is proved by a very small modification of Hardouin's arguments for \cite[Theorem 2]{Har11} and \cite[Theorem 2.1]{Har06}, those results in turn building on Bertrand's \cite[Theorem 1.1]{Be01}), we have
\[
D'\, ^p(\mu)_{A} \ = \ \{s\in S(\mu)(\CC)_{(x,y)} \, : \, \inHom((M/W_pM)_y, M_s)^\dagger/A \, \in \, \langle (W_pM)_x, (M/W_pM)_y\rangle^{\otimes} \}.
\]
Here, $\langle (W_pM)_x, (M/W_pM)_y\rangle^{\otimes}$ means the Tannakian subcategory of the category of mixed Hodge structures generated by $(W_pM)_x$ and $(M/W_pM)_y$. (By definition, this subcategory is closed under taking subobjects). 

Up to isomorphism, the category $\langle (W_pM)_x, (M/W_pM)_y\rangle ^{\otimes}$ has countably many objects (since for any finite-dimensional vector space $V$ and $H\leq GL(V)$, every representation of $H$ has a copy contained in a direct sum of the canonical representations $V^{\otimes m}\otimes {V\dual}^{\otimes n}$; see \cite[\S 4.14]{Mi17}, for instance). For any object $N$ of $\langle (W_pM)_x, (M/W_pM)_y\rangle^{\otimes}$, let
\[
D' \, ^{p}(\mu)_{A\, , \, N} \ := \ \{s\in S(\mu)(\CC)_{(x,y)} \, : \, \inHom((M/W_pM)_y, M_s)^\dagger/A \, \simeq \, N \}.
\] 
It is enough to show that each $D'\, ^p(\mu)_{A\, , \, N}$ is a union of countably many proper Zariski closed subsets of $S(\mu)(\CC)_{(x,y)} $. Note that since $A$ is a proper subobject of $\inHom((M/W_pM)_y, (W_pM)_x)$, the subset $D'\, ^p(\mu)_{A, \mathbbm{1}}$ is empty. Also note that $D'\, ^p(\mu)_{A\, , \, N}$ depends on the isomorphism class of $N$, rather than $N$ itself. 

We now fix an object $N\not\simeq \mathbbm{1}$ of $\langle (W_pM)_x, (M/W_pM)_y\rangle^{\otimes}$ and show that $D'\, ^p(\mu)_{A\, , \, N}$ is a union of countably many proper Zariski closed subsets of $S(\mu)(\CC)_{(x,y)}$. We shall assume that $D'\, ^p(\mu)_{A\, , \, N} $ is nonempty. There are only countably many filtrations $F^\cdot$
on 
\begin{equation}\label{eq1}
Hom(M_\CC/W_pM_\CC, M_\CC)^\dagger/A_\CC
\end{equation}
( = the underlying complex vector space of every $\inHom((M/W_pM)_y, M_s)^\dagger/A$) such that 
\[
(Hom(M_\QQ/W_pM_\QQ, M_\QQ)^\dagger/A_\QQ, W_\cdot, F^\cdot)
\]
(with $W_\cdot$ induced by our fixed weight filtration on $M_\QQ$) is a mixed Hodge structure isomorphic to $N$ (because after all, any such filtration will be the transport of the Hodge filtration on $N$ by some isomorphism $N_\QQ\longrightarrow Hom(M_\QQ/W_pM_\QQ, M_\QQ)^\dagger/A_\QQ$). Thus we may fix one such Hodge filtration $F_0^\cdot$ on Eq. \eqref{eq1} and show that
\begin{align*}
D' \, ^p(\mu)_{A\, , \, F_0^\cdot} \ := \ \bigm\{s\in S(\mu)(\CC)_{(x,y)} \, : \, & F^\cdot (\inHom((M/W_pM)_y, M_s)^\dagger/A) \\
\, & = \, F_0^\cdot \, (Hom(M_\CC/W_pM_\CC, M_\CC)^\dagger/A_\CC)\bigm\}
\end{align*}
is a proper Zariski closed subset of $S(\mu)(\CC)_{(x,y)}$. (We dropped $N$ from the notation as its isomorphism class is determined by $F_0^\cdot$\,.) 

\subsection{A construction} Before we prove that each $D' \, ^p(\mu)_{A\, , \, F_0^\cdot}$ is a proper Zariski closed subset of $S(\mu)(\CC)_{(x,y)}$ we need to give a construction. Let $\nu'$ be as in \S \ref{par: setup for proof of main thm}. Set 
\[
\wt{A} \ := \ Gr^WA \ \subsetneq \ \inHom(\wt{M}/W_p\wt{M}, W_p\wt{M}).
\] 
For any $s\in S(\mu)(\CC)$ such that $A_\QQ$ underlines a subobject of $\inHom((M/W_pM)_y, M_s)^\dagger$, and in particular for every $s\in S(\mu)(\CC)_{(x,y)}$, the mixed Hodge structure $\inHom((M/W_pM)_y, M_s)^\dagger/A$ is a mixed Hodge structure associated to the triple
\[
\nu'/{A_\QQ} \ := \ \bigm(Hom(M_\QQ/W_pM_\QQ, M_\QQ)^\dagger/A_\QQ, \, W_\cdot, \, \inHom(\wt{M}/W_p\wt{M}, \wt{M})^\dagger/ \wt{A}\bigm).
\]
The weight filtration of $Hom(M_\QQ/W_pM_\QQ, M_\QQ)^\dagger/A_\QQ$ has only weights $\leq 0$; moreover,
\[
W_{-1}Hom(M_\QQ/W_pM_\QQ, M_\QQ)^\dagger/A_\QQ \ = \ Hom(M_\QQ/W_pM_\QQ, W_pM_\QQ)/A_\QQ \ \neq \ 0
\]
and
\[
Gr^W_0(Hom(M_\QQ/W_pM_\QQ, M_\QQ)^\dagger/A_\QQ) \ = \ \QQ
\]
(via $\lambda$). With the notation of \S \ref{par: def of fibration}, we have 
\[
(\nu'/{A_\QQ})_{-1} \ = \ \bigm(Hom(M_\QQ/W_pM_\QQ, W_pM_\QQ)/A_\QQ, \, W_\cdot , \, \inHom(\wt{M}/W_p\wt{M}, W_p\wt{M})/ \wt{A}\bigm)
\]
and 
\[
(\nu'/{A_\QQ})_{>-1} \ = \ (\QQ, W_\cdot \mathbbm{1}, \mathbbm{1}),
\]
where $W_\cdot \mathbbm{1}$ means the weight filtration of the unit object $\mathbbm{1}$, considered at the rational level (i.e. the increasing filtration on $\QQ$ concentrated in step 0). Note that the only mixed Hodge structure associated to the latter triple is $\mathbbm{1}$ itself. Denote the element of $S((\nu'/{A_\QQ})_{-1})(\CC)$ corresponding to $\inHom((M/W_pM)_y, (W_pM)_x)/A$ by $(x;y)$. Let $S(\nu'/{A_\QQ})_{((x;y),\ast)}$ be the fiber of 
\[
\Theta_{-1} \, : \, S(\nu'/{A_\QQ}) \ \longrightarrow \ S((\nu'/{A_\QQ})_{-1}) \, \times \, S((\nu'/{A_\QQ})_{>-1}) \ = \ S((\nu'/{A_\QQ})_{-1}) \, \times \, \{\ast\}
\]
(i.e. $\Theta_{-1}$ for the triple $\nu'/{A_\QQ}$) above $((x;y),\ast)$. (See \S \ref{par: def of fibration}.)

We have a function
\[
\Delta \, : \, S(\mu)(\CC)_{(x,y)} \ \longrightarrow \ S(\nu'/{A_\QQ})(\CC)_{((x;y),\ast)} \hspace{.4in}s \ \mapsto \  \inHom((M/W_pM)_y, M_s)^\dagger/A,
\]
which in fact comes from a morphism
\[
S(\mu)_{(x,y)} \ \longrightarrow \ S(\nu'/{A_\QQ})_{((x;y),\ast)}.
\]
This can be easily seen using Proposition \ref{prop: fibers of theta_p}, but it can also be seen using the material of \S \ref{sec: algebraic families attached to triples}: Consider the morphism $S(\mu)\longrightarrow S(\nu')$ corresponding to the family 
\[(\inHom(M_s/W_pM_s, M_s)^\dagger)_{s\in S(\mu)(\CC)}\] 
(see Lemma \ref{alg families coming from linear algebra}); by Proposition \ref{prop quotients for alg families} it gives a morphism
\begin{equation}\label{eq39}
S(\mu)_{A_\QQ} \ \mapsto \ S(\nu'/A_\QQ),
\end{equation}
corresponding to the family $(\inHom(M_s/W_pM_s, M_s)^\dagger/A_s)_{s\in S(\mu)(\CC)_{A_\QQ}}$, where $S(\mu)_{A_\QQ}$ is with respect to the family $(\inHom(M_s/W_pM_s, M_s)^\dagger)_{s\in S(\mu)(\CC)}$. (By definition, it is the subvariety of $S(\mu)$ whose set of closed points $S(\mu)(\CC)_{A_\QQ}$ is the set of those $s\in S(\mu)(\CC)$ such that $A_\QQ$ underlies a subobject of $\inHom(M_s/W_pM_s, M_s)^\dagger$. By construction, $A_\QQ$ is a global Hodge subspace for the family $S(\mu)_{A_\QQ}\hookrightarrow S(\mu)\rightarrow S(\nu')$.) Since for every $s\in S(\mu)(\CC)_{(x,y)}$ we have
\begin{align*}
W_{-1}\inHom(M_s/W_pM_s, M_s)^\dagger/A_s \ &= \ \inHom(M_s/W_pM_s, W_pM_s)/A_s \\
 & = \ \inHom((M/W_pM)_y, (W_pM)_x)/A,
\end{align*}
the morphism Eq. \eqref{eq39} restricts to a morphism
\[
S(\mu)_{(x,y)} \ \longrightarrow \ S(\nu'/{A_\QQ})_{((x;y),\ast)}.
\]
At the complex level this is the map $\Delta$.

\subsection{Finishing the proof} 
We now come back to showing that $D' \, ^p(\mu)_{A\, , \, F_0^\cdot}$ is a proper Zariski closed subset of $S(\mu)(\CC)_{(x,y)}$ (with $A$ is as before, and $F_0$ as in \S \ref{par: setup proof of Theorem C}: a fixed filtration on Eq. \eqref{eq1} making $Hom(M_\QQ/W_pM_\QQ, M_\QQ)^\dagger/A_\QQ$ with its given weight filtration a mixed Hodge structure). We will assume that $D' \, ^p(\mu)_{A\, , \, F_0^\cdot}$ is nonempty. Since $F^\cdot_0$ is the Hodge filtration on $\inHom(M_s/W_pM_s, M_s)^\dagger/A_s$ for some $s\in S(\mu)(\CC)_{(x,y)}$, it follows that the mixed Hodge structure on $Hom(M_\QQ/W_pM_\QQ, M_\QQ)^\dagger/A_\QQ$ with Hodge filtration $F^\cdot_0$ (and the given weight filtration) is associated to the triple $\nu'/A_\QQ$, and moreover belongs to $S(\nu'/{A_\QQ})(\CC)_{((x;y),\ast)}$. Now the set $D' \, ^p(\mu)_{A\, , \, F_0^\cdot}$ is exactly the pre-image of this element of $S(\nu'/{A_\QQ})(\CC)_{((x;y),\ast)}$ under the map $\Delta$, and hence is Zariski closed in $S(\mu)(\CC)_{(x,y)}$.

It remains to establish properness. This is done by verifying that (i) $\Delta$ is surjective and (ii) $S(\nu'/{A_\QQ})(\CC)_{((x;y),\ast)}$ is more than a single point. By Proposition \ref{prop: fibers of theta_p} (applied to $\Theta_{-1}$ of $\nu'/A_\QQ$), every element of $S(\nu'/{A_\QQ})(\CC)_{((x;y),\ast)}$ is of the form $[\overline{\psi}]$ for a section $\overline{\psi}$ of the map
\[
\overline{\lambda} \, : \, Hom(M_\CC/W_pM_\CC,M_\CC)^\dagger/A_\CC \ \longrightarrow \ \CC.
\]
Here, $\overline{\lambda}$ is the map induced by $\lambda: Hom(M_\CC/W_pM_\CC,M_\CC)^\dagger \longrightarrow  \CC$, and $[\overline{\psi}]$ denotes the element of the fiber $S(\nu'/{A_\QQ})(\CC)_{((x;y),\ast)}$ arising from the section $\overline{\psi}$ via the isomorphism of Proposition \ref{prop: fibers of theta_p}. Now lift $\overline{\psi}$ to section $\psi$ of $\lambda$. Then $\psi(1) \in Hom(M_\CC/W_pM_\CC,M_\CC)^\dagger$ is a section of $M_\CC\longrightarrow M_\CC/W_pM_\CC$, and hence by Proposition \ref{prop: fibers of theta_p} (applied to $\Theta_p$ for $\mu$) gives an element $[\psi(1)]$ of $S(\mu)(\CC)_{(x,y)}$. Then $\Delta$ sends $[\psi(1)]$ to $[\overline{\psi}]$.

Finally, note that again by Proposition \ref{prop: fibers of theta_p},
\begin{align*}
S(\nu'/{A_\QQ})(\CC)_{((x;y),\ast)} \  \cong \ \ \ \ \ \ \  & \ \frac{\overline{\lambda}^{-1}(1)}{F^0(\inHom((M/W_pM)_y,(W_pM)_x)/A)}\\
 \stackrel{\text{non-canonically}}{\simeq} & \ \frac{Hom(M_\CC/W_pM_\CC, W_pM_\CC)/A_\CC}{F^0(\inHom((M/W_pM)_y,(W_pM)_x)/A)}.
\end{align*}
The object  
\[
\inHom((M/W_pM)_y, (W_pM)_x)/A
\]
is a nonzero object with all weights $<0$, and hence its $F^0$ cannot be the entire space. This completes the proof of Theorem \ref{thm1}.

\end{document}